%
%
%
%
%
%
\year=1996
%
%
%
%
\expandafter\ifx\csname TeX\endcsname\relax \input plain \fi
\expandafter\ifx\csname phyzzx\endcsname\relax \else
 \wlog{PHYZZX macros are already loaded and are not
          \string\input\space again}%
 \endinput \fi
\catcode`\@=11 
\let\rel@x=\relax
\let\n@expand=\relax
\def\pr@tect{\let\n@expand=\noexpand}
\let\protect=\pr@tect
\let\gl@bal=\global
%
%
%
\newfam\cpfam
\newdimen\b@gheight             \b@gheight=12pt
\newcount\f@ntkey               \f@ntkey=0
\def\f@m{\afterassignment\samef@nt\f@ntkey=}
\def\samef@nt{\fam=\f@ntkey\the\textfont\f@ntkey\rel@x}
\def\setstr@t{\setbox\strutbox=\hbox{\vrule height 0.85\b@gheight
                                depth 0.35\b@gheight width\z@ }}
%
%
%
%

\font\fourteenrm  =cmr12 scaled\magstep1
\font\twelverm    =cmr12
\font\ninerm      =cmr9
\font\sixrm       =cmr6

\font\fourteenbf  =cmbx12 scaled\magstep1
\font\twelvebf    =cmbx12
\font\ninebf      =cmbx9
\font\sixbf       =cmbx6
\font\seventeeni  =cmmi12 scaled\magstep2    \skewchar\seventeeni='177
\font\fourteeni   =cmmi12 scaled\magstep1     \skewchar\fourteeni='177
\font\twelvei     =cmmi12                       \skewchar\twelvei='177
\font\ninei       =cmmi9                          \skewchar\ninei='177
\font\sixi        =cmmi6                           \skewchar\sixi='177
\font\seventeensy =cmsy10 scaled\magstep3    \skewchar\seventeensy='60
\font\fourteensy  =cmsy10 scaled\magstep2     \skewchar\fourteensy='60
\font\twelvesy    =cmsy10 scaled\magstep1       \skewchar\twelvesy='60
\font\ninesy      =cmsy9                          \skewchar\ninesy='60
\font\sixsy       =cmsy6                           \skewchar\sixsy='60

\font\fourteenex  =cmex10 scaled\magstep2
\font\twelveex    =cmex10 scaled\magstep1

\font\fourteensl  =cmsl12 scaled\magstep1
\font\twelvesl    =cmsl12
\font\ninesl      =cmsl9

\font\fourteenit  =cmti12 scaled\magstep1
\font\twelveit    =cmti12
\font\nineit      =cmti9
\font\fourteentt  =cmtt12 scaled\magstep1
\font\twelvett    =cmtt12
\font\fourteencp  =cmcsc10 scaled\magstep2
\font\twelvecp    =cmcsc10 scaled\magstep1
\font\tencp       =cmcsc10
%
%
\def\fourteenf@nts{\relax
    \textfont0=\fourteenrm          \scriptfont0=\tenrm
      \scriptscriptfont0=\sevenrm
    \textfont1=\fourteeni           \scriptfont1=\teni
      \scriptscriptfont1=\seveni
    \textfont2=\fourteensy          \scriptfont2=\tensy
      \scriptscriptfont2=\sevensy
    \textfont3=\fourteenex          \scriptfont3=\twelveex
      \scriptscriptfont3=\tenex
    \textfont\itfam=\fourteenit     \scriptfont\itfam=\tenit
    \textfont\slfam=\fourteensl     \scriptfont\slfam=\tensl
    \textfont\bffam=\fourteenbf     \scriptfont\bffam=\tenbf
      \scriptscriptfont\bffam=\sevenbf
    \textfont\ttfam=\fourteentt
    \textfont\cpfam=\fourteencp }
\def\twelvef@nts{\relax
    \textfont0=\twelverm          \scriptfont0=\ninerm
      \scriptscriptfont0=\sixrm
    \textfont1=\twelvei           \scriptfont1=\ninei
      \scriptscriptfont1=\sixi
    \textfont2=\twelvesy           \scriptfont2=\ninesy
      \scriptscriptfont2=\sixsy
    \textfont3=\twelveex          \scriptfont3=\tenex
      \scriptscriptfont3=\tenex
    \textfont\itfam=\twelveit     \scriptfont\itfam=\nineit
    \textfont\slfam=\twelvesl     \scriptfont\slfam=\ninesl
    \textfont\bffam=\twelvebf     \scriptfont\bffam=\ninebf
      \scriptscriptfont\bffam=\sixbf
    \textfont\ttfam=\twelvett
    \textfont\cpfam=\twelvecp }
\def\tenf@nts{\relax
    \textfont0=\tenrm          \scriptfont0=\sevenrm
      \scriptscriptfont0=\fiverm
    \textfont1=\teni           \scriptfont1=\seveni
      \scriptscriptfont1=\fivei
    \textfont2=\tensy          \scriptfont2=\sevensy
      \scriptscriptfont2=\fivesy
    \textfont3=\tenex          \scriptfont3=\tenex
      \scriptscriptfont3=\tenex
    \textfont\itfam=\tenit     \scriptfont\itfam=\seveni  
    \textfont\slfam=\tensl     \scriptfont\slfam=\sevenrm 
    \textfont\bffam=\tenbf     \scriptfont\bffam=\sevenbf
      \scriptscriptfont\bffam=\fivebf
    \textfont\ttfam=\tentt
    \textfont\cpfam=\tencp }
%
%
%
%
\def\rm{\n@expand\f@m0 }
\def\mit{\n@expand\f@m1 }         \let\oldstyle=\mit
\def\cal{\n@expand\f@m2}
\def\it{\n@expand\f@m\itfam}
\def\sl{\n@expand\f@m\slfam}
\def\bf{\n@expand\f@m\bffam}
\def\tt{\n@expand\f@m\ttfam}
\def\caps{\n@expand\f@m\cpfam}    \let\cp=\caps
\def\em@{\rel@x\ifnum\f@ntkey=0\it\else
        \ifnum\f@ntkey=\bffam\it\else\rm\fi \fi }
\def\em{\n@expand\em@}
\def\fourteenpoint{\fourteenf@nts \samef@nt \b@gheight=14pt \setstr@t }
\def\twelvepoint{\twelvef@nts \samef@nt \b@gheight=12pt \setstr@t }
\def\tenpoint{\tenf@nts \samef@nt \b@gheight=10pt \setstr@t }
\normalbaselineskip = 20pt plus 0.2pt minus 0.1pt
\normallineskip = 1.5pt plus 0.1pt minus 0.1pt
\normallineskiplimit = 1.5pt
\newskip\normaldisplayskip
\normaldisplayskip = 20pt plus 5pt minus 10pt
\newskip\normaldispshortskip
\normaldispshortskip = 6pt plus 5pt
\newskip\normalparskip
\normalparskip = 6pt plus 2pt minus 1pt
\newskip\skipregister
\skipregister = 5pt plus 2pt minus 1.5pt
\newif\ifsingl@
\newif\ifdoubl@
\newif\iftwelv@  \twelv@true
\def\singlespace{\singl@true\doubl@false\spaces@t}
\def\doublespace{\singl@false\doubl@true\spaces@t}
\def\normalspace{\singl@false\doubl@false\spaces@t}
\def\Tenpoint{\tenpoint\twelv@false\spaces@t}
\def\Twelvepoint{\twelvepoint\twelv@true\spaces@t}
\def\spaces@t{\rel@x
      \iftwelv@ \ifsingl@\subspaces@t3:4;\else\subspaces@t1:1;\fi
       \else \ifsingl@\subspaces@t3:5;\else\subspaces@t4:5;\fi \fi
      \ifdoubl@ \multiply\baselineskip by 5
         \divide\baselineskip by 4 \fi }
\def\subspaces@t#1:#2;{
      \baselineskip = \normalbaselineskip
      \multiply\baselineskip by #1 \divide\baselineskip by #2
      \lineskip = \normallineskip
      \multiply\lineskip by #1 \divide\lineskip by #2
      \lineskiplimit = \normallineskiplimit
      \multiply\lineskiplimit by #1 \divide\lineskiplimit by #2
      \parskip = \normalparskip
      \multiply\parskip by #1 \divide\parskip by #2
      \abovedisplayskip = \normaldisplayskip
      \multiply\abovedisplayskip by #1 \divide\abovedisplayskip by #2
      \belowdisplayskip = \abovedisplayskip
      \abovedisplayshortskip = \normaldispshortskip
      \multiply\abovedisplayshortskip by #1
        \divide\abovedisplayshortskip by #2
      \belowdisplayshortskip = \abovedisplayshortskip
      \advance\belowdisplayshortskip by \belowdisplayskip
      \divide\belowdisplayshortskip by 2
      \smallskipamount = \skipregister
      \multiply\smallskipamount by #1 \divide\smallskipamount by #2
      \medskipamount = \smallskipamount \multiply\medskipamount by 2
      \bigskipamount = \smallskipamount \multiply\bigskipamount by 4 }
\def\normalbaselines{ \baselineskip=\normalbaselineskip
   \lineskip=\normallineskip \lineskiplimit=\normallineskip
   \iftwelv@\else \multiply\baselineskip by 4 \divide\baselineskip by 5
     \multiply\lineskiplimit by 4 \divide\lineskiplimit by 5
     \multiply\lineskip by 4 \divide\lineskip by 5 \fi }
\Twelvepoint  
\interlinepenalty=50
\interfootnotelinepenalty=5000
\predisplaypenalty=9000
\postdisplaypenalty=500
\hfuzz=1pt
\vfuzz=0.2pt
\newdimen\HOFFSET  \HOFFSET=0pt
\newdimen\VOFFSET  \VOFFSET=0pt
\newdimen\HSWING   \HSWING=0pt
\dimen\footins=8in
%
%
%
\newskip\pagebottomfiller
\pagebottomfiller=\z@ plus \z@ minus \z@
\def\pagecontents{
   \ifvoid\topins\else\unvbox\topins\vskip\skip\topins\fi
   \dimen@ = \dp255 \unvbox255
   \vskip\pagebottomfiller
   \ifvoid\footins\else\vskip\skip\footins\footrule\unvbox\footins\fi
   \ifr@ggedbottom \kern-\dimen@ \vfil \fi }
\def\makeheadline{\vbox to 0pt{ \skip@=\topskip
      \advance\skip@ by -12pt \advance\skip@ by -2\normalbaselineskip
      \vskip\skip@ \line{\vbox to 12pt{}\the\headline} \vss
      }\nointerlineskip}
\def\makefootline{\baselineskip = 1.5\normalbaselineskip
                 \line{\the\footline}}
\newif\iffrontpage
\newif\ifp@genum
\def\nopagenumbers{\p@genumfalse}
\def\pagenumbers{\p@genumtrue}
\pagenumbers
\newtoks\paperheadline
\newtoks\paperfootline
\newtoks\letterheadline
\newtoks\letterfootline
\newtoks\letterinfo
\newtoks\date
\paperheadline={\hfil}
\paperfootline={\hss\iffrontpage\else\ifp@genum\tenrm\folio\hss\fi\fi}
\letterheadline{\iffrontpage \hfil \else
    \rm \ifp@genum page~~\folio\fi \hfil\the\date \fi}
\letterfootline={\iffrontpage\the\letterinfo\else\hfil\fi}
\letterinfo={\hfil}
\def\monthname{\rel@x\ifcase\month 0/\or January\or February\or
   March\or April\or May\or June\or July\or August\or September\or
   October\or November\or December\else\number\month/\fi}
\def\today{\monthname~\number\day, \number\year}
\date={\today}
\headline=\paperheadline 
\footline=\paperfootline 
\countdef\pageno=1      \countdef\pagen@=0
\countdef\pagenumber=1  \pagenumber=1
\def\advancepageno{\gl@bal\advance\pagen@ by 1
   \ifnum\pagenumber<0 \gl@bal\advance\pagenumber by -1
    \else\gl@bal\advance\pagenumber by 1 \fi
    \gl@bal\frontpagefalse  \swing@ }
\def\folio{\ifnum\pagenumber<0 \romannumeral-\pagenumber
           \else \number\pagenumber \fi }
\def\swing@{\ifodd\pagenumber \gl@bal\advance\hoffset by -\HSWING
             \else \gl@bal\advance\hoffset by \HSWING \fi }
\def\footrule{\dimen@=\prevdepth\nointerlineskip
   \vbox to 0pt{\vskip -0.25\baselineskip \hrule width 0.35\hsize \vss}
   \prevdepth=\dimen@ }
\let\footnotespecial=\rel@x
\newdimen\footindent
\footindent=24pt
\def\Textindent#1{\noindent\llap{#1\enspace}\ignorespaces}
\def\Vfootnote#1{\insert\footins\bgroup
   \interlinepenalty=\interfootnotelinepenalty \floatingpenalty=20000
   \singl@true\doubl@false\Tenpoint
   \splittopskip=\ht\strutbox \boxmaxdepth=\dp\strutbox
   \leftskip=\footindent \rightskip=\z@skip
   \parindent=0.5\footindent \parfillskip=0pt plus 1fil
   \spaceskip=\z@skip \xspaceskip=\z@skip \footnotespecial
   \Textindent{#1}\footstrut\futurelet\next\fo@t}

\def\vfootnote#1{\Vfootnote{${#1}$}}
\def\footnote#1{\attach{#1}\vfootnote{#1}}

\def\foot{\attach\footsymbolgen\vfootnote{\footsymbol}}
\let\footsymbol=\star
\newcount\lastf@@t           \lastf@@t=-1
\newcount\footsymbolcount    \footsymbolcount=0
\newif\ifPhysRev
\def\footsymbolgen{\bumpfootsymbolcount \generatefootsymbol \footsymbol }
\def\bumpfootsymbolcount{\rel@x
   \iffrontpage \bumpfootsymbolpos \else \advance\lastf@@t by 1
     \ifPhysRev \bumpfootsymbolneg \else \bumpfootsymbolpos \fi \fi
   \gl@bal\lastf@@t=\pagen@ }
\def\bumpfootsymbolpos{\ifnum\footsymbolcount <0
                            \gl@bal\footsymbolcount =0 \fi
    \ifnum\lastf@@t<\pagen@ \gl@bal\footsymbolcount=0
     \else \gl@bal\advance\footsymbolcount by 1 \fi }
\def\bumpfootsymbolneg{\ifnum\footsymbolcount >0
             \gl@bal\footsymbolcount =0 \fi
         \gl@bal\advance\footsymbolcount by -1 }
\def\fd@f#1 {\xdef\footsymbol{\mathchar"#1 }}
\def\generatefootsymbol{\ifcase\footsymbolcount \fd@f 13F \or \fd@f 279
        \or \fd@f 27A \or \fd@f 278 \or \fd@f 27B \else
        \ifnum\footsymbolcount <0 \fd@f{023 \number-\footsymbolcount }
         \else \fd@f 203 {\loop \ifnum\footsymbolcount >5
                \fd@f{203 \footsymbol } \advance\footsymbolcount by -1
                \repeat }\fi \fi }

\def\nonfrenchspacing{\sfcode`\.=3001 \sfcode`\!=3000 \sfcode`\?=3000
        \sfcode`\:=2000 \sfcode`\;=1500 \sfcode`\,=1251 }
\nonfrenchspacing
\newdimen\d@twidth
{\setbox0=\hbox{s.} \gl@bal\d@twidth=\wd0 \setbox0=\hbox{s}
        \gl@bal\advance\d@twidth by -\wd0 }
\def\removehglue{\loop \unskip \ifdim\lastskip >\z@ \repeat }
\def\roll@ver#1{\removehglue \nobreak \count255 =\spacefactor \dimen@=\z@
        \ifnum\count255 =3001 \dimen@=\d@twidth \fi
        \ifnum\count255 =1251 \dimen@=\d@twidth \fi
    \iftwelv@ \kern-\dimen@ \else \kern-0.83\dimen@ \fi
   #1\spacefactor=\count255 }
\def\step@ver#1{\rel@x \ifmmode #1\else \ifhmode
        \roll@ver{${}#1$}\else {\setbox0=\hbox{${}#1$}}\fi\fi }
\def\attach#1{\step@ver{\strut^{\mkern 2mu #1} }}
%
%
%
\newcount\chapternumber      \chapternumber=0
\newcount\sectionnumber      \sectionnumber=0
\newcount\equanumber         \equanumber=0
\let\chapterlabel=\rel@x
\let\sectionlabel=\rel@x
\newtoks\chapterstyle        \chapterstyle={\Number}
\newtoks\sectionstyle        \sectionstyle={\Number}
\newskip\chapterskip         \chapterskip=\bigskipamount
\newskip\sectionskip         \sectionskip=\medskipamount
\newskip\headskip            \headskip=8pt plus 3pt minus 3pt
\newdimen\chapterminspace    \chapterminspace=15pc
\newdimen\sectionminspace    \sectionminspace=10pc
\newdimen\referenceminspace  \referenceminspace=20pc
\newif\ifcn@                 \cn@true
\newif\ifcn@@                \cn@@false
\def\numberedchapters{\cn@true}
\def\unnumberedchapters{\cn@false\sequentialequations}
\def\chapterreset{\gl@bal\advance\chapternumber by 1
   \ifnum\equanumber<0 \else\gl@bal\equanumber=0\fi
   \sectionnumber=0 \let\sectionlabel=\rel@x
   \ifcn@ \gl@bal\cn@@true {\pr@tect
       \xdef\chapterlabel{\the\chapterstyle{\the\chapternumber}}}%
    \else \gl@bal\cn@@false \gdef\chapterlabel{\rel@x}\fi }
\def\@alpha#1{\count255='140 \advance\count255 by #1\char\count255}
 \def\alphabetic{\n@expand\@alpha}
\def\@Alpha#1{\count255='100 \advance\count255 by #1\char\count255}
 \def\Alphabetic{\n@expand\@Alpha}
\def\@Roman#1{\uppercase\expandafter{\romannumeral #1}}
 \def\Roman{\n@expand\@Roman}
\def\@roman#1{\romannumeral #1}    \def\roman{\n@expand\@roman}
\def\@number#1{\number #1}         \def\Number{\n@expand\@number}
\def\BLANK#1{\rel@x}               
\def\titleparagraphs{\interlinepenalty=9999
     \leftskip=0.03\hsize plus 0.22\hsize minus 0.03\hsize
     \rightskip=\leftskip \parfillskip=0pt
     \hyphenpenalty=9000 \exhyphenpenalty=9000
     \tolerance=9999 \pretolerance=9000
     \spaceskip=0.333em \xspaceskip=0.5em }
\def\titlestyle#1{\par\begingroup \titleparagraphs
     \iftwelv@\fourteenpoint\else\twelvepoint\fi
   \noindent #1\par\endgroup }
\def\spacecheck#1{\dimen@=\pagegoal\advance\dimen@ by -\pagetotal
   \ifdim\dimen@<#1 \ifdim\dimen@>0pt \vfil\break \fi\fi}
\def\chapter#1{\par \penalty-300 \vskip\chapterskip
   \spacecheck\chapterminspace
   \chapterreset \titlestyle{\ifcn@@\chapterlabel.~\fi #1}
   \nobreak\vskip\headskip \penalty 30000
   {\pr@tect\wlog{\string\chapter\space \chapterlabel}} }

\def\section#1{\par \ifnum\lastpenalty=30000\else
   \penalty-200\vskip\sectionskip \spacecheck\sectionminspace\fi
   \gl@bal\advance\sectionnumber by 1
   {\pr@tect
   \xdef\sectionlabel{\ifcn@@ \chapterlabel.\fi
       \the\sectionstyle{\the\sectionnumber}}%
   \wlog{\string\section\space \sectionlabel}}%
   \noindent {\caps\enspace\sectionlabel.~~#1}\par
   \nobreak\vskip\headskip \penalty 30000 }
\def\subsection#1{\par
   \ifnum\the\lastpenalty=30000\else \penalty-100\smallskip \fi
   \noindent\undertext{#1}\enspace \vadjust{\penalty5000}}

\def\undertext#1{\vtop{\hbox{#1}\kern 1pt \hrule}}

\def\ack{\subsection{Acknowledgements:}}
\def\APPENDIX#1#2{\par\penalty-300\vskip\chapterskip
   \spacecheck\chapterminspace \chapterreset \xdef\chapterlabel{#1}
   \titlestyle{APPENDIX #2} \nobreak\vskip\headskip \penalty 30000
   \wlog{\string\Appendix~\chapterlabel} }
\def\Appendix#1{\APPENDIX{#1}{#1}}
\def\appendix{\APPENDIX{A}{}}
%
%
%
%

\def\eqn{\eqno\eqname}

\def\eqinsert#1{\noalign{\dimen@=\prevdepth \nointerlineskip
   \setbox0=\hbox to\displaywidth{\hfil #1}
   \vbox to 0pt{\kern 0.5\baselineskip\hbox{$\!\box0\!$}\vss}
   \prevdepth=\dimen@}}
%

%
%
\def\GENITEM#1;#2{\par \hangafter=0 \hangindent=#1
    \Textindent{$ #2 $}\ignorespaces}
\outer\def\newitem#1=#2;{\gdef#1{\GENITEM #2;}}

\newdimen\itemsize                \itemsize=30pt
\newitem\item=1\itemsize;
\newitem\sitem=1.75\itemsize;     
\newitem\ssitem=2.5\itemsize;     
\outer\def\newlist#1=#2&#3&#4;{\toks0={#2}\toks1={#3}%
   \count255=\escapechar \escapechar=-1
   \alloc@0\list\countdef\insc@unt\listcount     \listcount=0
   \edef#1{\par
      \countdef\listcount=\the\allocationnumber
      \advance\listcount by 1
      \hangafter=0 \hangindent=#4
      \Textindent{\the\toks0{\listcount}\the\toks1}}
   \expandafter\expandafter\expandafter
    \edef\c@t#1{begin}{\par
      \countdef\listcount=\the\allocationnumber \listcount=1
      \hangafter=0 \hangindent=#4
      \Textindent{\the\toks0{\listcount}\the\toks1}}
   \expandafter\expandafter\expandafter
    \edef\c@t#1{con}{\par \hangafter=0 \hangindent=#4 \noindent}
   \escapechar=\count255}
\def\c@t#1#2{\csname\string#1#2\endcsname}
\newlist\point=\Number&.&1.0\itemsize;
\newlist\subpoint=(\alphabetic&)&1.75\itemsize;
\newlist\subsubpoint=(\roman&)&2.5\itemsize;
%

%
%
%
%
\newcount\referencecount     \referencecount=0
\newcount\lastrefsbegincount \lastrefsbegincount=0
\newif\ifreferenceopen       \newwrite\referencewrite
\newdimen\refindent          \refindent=30pt
\def\normalrefmark#1{\attach{\scriptscriptstyle [ #1 ] }}
\let\PRrefmark=\attach
\def\NPrefmark#1{\step@ver{{\;[#1]}}}
\def\refmark#1{\rel@x\ifPhysRev\PRrefmark{#1}\else\normalrefmark{#1}\fi}
\def\refend@{\refmark{\number\referencecount}}
\def\refend{\refend@{}\space }
\def\refsend{\refmark{\count255=\referencecount
   \advance\count255 by-\lastrefsbegincount
   \ifcase\count255 \number\referencecount
   \or \number\lastrefsbegincount,\number\referencecount
   \else \number\lastrefsbegincount-\number\referencecount \fi}\space }
\def\REFNUM#1{\rel@x \gl@bal\advance\referencecount by 1
    \xdef#1{\the\referencecount }}
\def\Refnum#1{\REFNUM #1\refend@ } 
\def\REF#1{\REFNUM #1\R@FWRITE\ignorespaces}
\def\Ref#1{\Refnum #1\REFWRITE }
\def\ref{\Ref\?}
\def\REFS#1{\REFNUM #1\gl@bal\lastrefsbegincount=\referencecount
    \REFWRITE }

\def\r@fitem#1{\par \hangafter=0 \hangindent=\refindent \Textindent{#1}}
\def\refitem#1{\r@fitem{#1.}}
\def\NPrefitem#1{\r@fitem{[#1]}}
\def\NPrefs{\let\refmark=\NPrefmark \let\refitem=NPrefitem}
\def\REFWRITE{\R@FWRITE\rel@x }
\def\R@FWRITE#1{\ifreferenceopen \else \gl@bal\referenceopentrue
     \immediate\openout\referencewrite=\jobname.refs
     \toks@={\begingroup \refoutspecials \catcode`\^^M=10 }%
     \immediate\write\referencewrite{\the\toks@}\fi
    \immediate\write\referencewrite{\noexpand\refitem %
                                    {\the\referencecount}}%
    \p@rse@ndwrite \referencewrite #1}
\begingroup
 \catcode`\^^M=\active \let^^M=\relax %
 \gdef\p@rse@ndwrite#1#2{\begingroup \catcode`\^^M=12 \newlinechar=`\^^M%
         \chardef\rw@write=#1\sc@nlines#2}%
 \gdef\sc@nlines#1#2{\sc@n@line \g@rbage #2^^M\endsc@n \endgroup #1}%
 \gdef\sc@n@line#1^^M{\expandafter\toks@\expandafter{\deg@rbage #1}%
         \immediate\write\rw@write{\the\toks@}%
         \futurelet\n@xt \sc@ntest }%
\endgroup
\def\sc@ntest{\ifx\n@xt\endsc@n \let\n@xt=\rel@x
       \else \let\n@xt=\sc@n@notherline \fi \n@xt }
\def\sc@n@notherline{\sc@n@line \g@rbage }
\def\deg@rbage#1{}
\let\g@rbage=\relax    \let\endsc@n=\relax
\def\refout{\par\penalty-400\vskip\chapterskip
   \spacecheck\referenceminspace
   \ifreferenceopen \Closeout\referencewrite \referenceopenfalse \fi
   \line{\fourteenrm\hfil REFERENCES\hfil}\vskip\headskip
   \input \jobname.refs
   }
\def\refoutspecials{\sfcode`\.=1000 \interlinepenalty=1000
         \rightskip=\z@ plus 1em minus \z@ }
\def\Closeout#1{\toks0={\par\endgroup}\immediate\write#1{\the\toks0}%
   \immediate\closeout#1}
%
%
\newcount\figurecount     \figurecount=0
\newcount\tablecount      \tablecount=0
\newif\iffigureopen       \newwrite\figurewrite
\newif\iftableopen        \newwrite\tablewrite
\def\FIGNUM#1{\rel@x \gl@bal\advance\figurecount by 1
    \xdef#1{\the\figurecount}}
\def\FIGURE#1{\FIGNUM #1\F@GWRITE\ignorespaces }

\def\figitem#1{\r@fitem{#1)}}
\def\FIGWRITE{\F@GWRITE\rel@x }
\def\TABNUM#1{\rel@x \gl@bal\advance\tablecount by 1
    \xdef#1{\the\tablecount}}
\def\TABLE#1{\TABNUM #1\T@BWRITE\ignorespaces }
\def\Table{\TABNUM\?Table~\?\TABWRITE }
\def\tabitem#1{\r@fitem{#1:}}
\def\TABWRITE{\T@BWRITE\rel@x }
\def\F@GWRITE#1{\iffigureopen \else \gl@bal\figureopentrue
     \immediate\openout\figurewrite=\jobname.figs
     \toks@={\begingroup \catcode`\^^M=10 }%
     \immediate\write\figurewrite{\the\toks@}\fi
    \immediate\write\figurewrite{\noexpand\figitem %
                                 {\the\figurecount}}%
    \p@rse@ndwrite \figurewrite #1}
\def\T@BWRITE#1{\iftableopen \else \gl@bal\tableopentrue
     \immediate\openout\tablewrite=\jobname.tabs
     \toks@={\begingroup \catcode`\^^M=10 }%
     \immediate\write\tablewrite{\the\toks@}\fi
    \immediate\write\tablewrite{\noexpand\tabitem %
                                 {\the\tablecount}}%
    \p@rse@ndwrite \tablewrite #1}
\def\figout{\par\penalty-400
   \vskip\chapterskip\spacecheck\referenceminspace
   \iffigureopen \Closeout\figurewrite \figureopenfalse \fi
   \line{\fourteenrm\hfil FIGURE CAPTIONS\hfil}\vskip\headskip
   \input \jobname.figs
   }
\def\tabout{\par\penalty-400
   \vskip\chapterskip\spacecheck\referenceminspace
   \iftableopen \Closeout\tablewrite \tableopenfalse \fi
   \line{\fourteenrm\hfil TABLE CAPTIONS\hfil}\vskip\headskip
   \input \jobname.tabs
   }
%
%
%
\newbox\picturebox
\def\p@cht{\ht\picturebox }
\def\p@cwd{\wd\picturebox }
\def\p@cdp{\dp\picturebox }
\newdimen\xshift
\newdimen\yshift
\newdimen\captionwidth
\newskip\captionskip
\captionskip=15pt plus 5pt minus 3pt
\def\fullwidth{\captionwidth=\hsize }
\newtoks\Caption
\newif\ifcaptioned
\newif\ifselfcaptioned
\def\caption{\captionedtrue \Caption }
\newcount\linesabove
\newif\iffileexists
\newtoks\picfilename
\def\fil@#1 {\fileexiststrue \picfilename={#1}}
\def\file#1{\if=#1\let\n@xt=\fil@ \else \def\n@xt{\fil@ #1}\fi \n@xt }
\def\pl@t{\begingroup \pr@tect
    \setbox\picturebox=\hbox{}\fileexistsfalse
    \let\height=\p@cht \let\width=\p@cwd \let\depth=\p@cdp
    \xshift=\z@ \yshift=\z@ \captionwidth=\z@
    \Caption={}\captionedfalse
    \linesabove =0 \picturedefault }
\def\plot{\pl@t \selfcaptionedfalse }
\def\Picture#1{\gl@bal\advance\figurecount by 1
    \xdef#1{\the\figurecount}\pl@t \selfcaptionedtrue }

\def\s@vepicture{\iffileexists \parsefilename \redopicturebox \fi
   \ifdim\captionwidth>\z@ \else \captionwidth=\p@cwd \fi
   \xdef\lastpicture{%
      \iffileexists%
         \setbox0=\hbox{\raise\the\yshift \vbox{%
              \moveright\the\xshift\hbox{\picturedefinition}}}%
      \else%
         \setbox0=\hbox{}%
      \fi
      \ht0=\the\p@cht \wd0=\the\p@cwd \dp0=\the\p@cdp
      \vbox{\hsize=\the\captionwidth%
            \line{\hss\box0 \hss }%
            \ifcaptioned%
               \vskip\the\captionskip \noexpand\Tenpoint
               \ifselfcaptioned%
                   Figure~\the\figurecount.\enspace%
               \fi%
               \the\Caption%
           \fi%
           }%
      }%
      \endgroup%
}
\let\endpicture=\s@vepicture
\def\savepicture#1{\s@vepicture \global\let#1=\lastpicture }
\def\displaypicture{\fullwidth \s@vepicture $$\lastpicture $${}}
\def\toppicture{\fullwidth \s@vepicture \topinsert
    \lastpicture \medskip \endinsert }
\def\midpicture{\fullwidth \s@vepicture \midinsert
    \lastpicture \endinsert }
%
%
\def\leftpicture{\pres@tpicture
    \dimen@i=\hsize \advance\dimen@i by -\dimen@ii
    \setbox\picturebox=\hbox to \hsize {\box0 \hss }%
    \wr@paround }
\def\rightpicture{\pres@tpicture
    \dimen@i=\z@
    \setbox\picturebox=\hbox to \hsize {\hss \box0 }%
    \wr@paround }
\def\pres@tpicture{\gl@bal\linesabove=\linesabove
    \s@vepicture \setbox\picturebox=\vbox{
         \kern \linesabove\baselineskip \kern 0.3\baselineskip
         \lastpicture \kern 0.3\baselineskip }%
    \dimen@=\p@cht \dimen@i=\dimen@
    \advance\dimen@i by \pagetotal
    \par \ifdim\dimen@i>\pagegoal \vfil\break \fi
    \dimen@ii=\hsize
    \advance\dimen@ii by -\parindent \advance\dimen@ii by -\p@cwd
    \setbox0=\vbox to\z@{\kern-\baselineskip \unvbox\picturebox \vss }}
\def\wr@paround{\Caption={}\count255=1
    \loop \ifnum \linesabove >0
         \advance\linesabove by -1 \advance\count255 by 1
         \advance\dimen@ by -\baselineskip
         \expandafter\Caption \expandafter{\the\Caption \z@ \hsize }%
      \repeat
    \loop \ifdim \dimen@ >\z@
         \advance\count255 by 1 \advance\dimen@ by -\baselineskip
         \expandafter\Caption \expandafter{%
             \the\Caption \dimen@i \dimen@ii }%
      \repeat
    \edef\n@xt{\parshape=\the\count255 \the\Caption \z@ \hsize }%
    \par\noindent \n@xt \strut \vadjust{\box\picturebox }}
\let\picturedefault=\relax
\let\parsefilename=\relax
\def\redopicturebox{\let\picturedefinition=\rel@x
   \errhelp=\disabledpictures
   \errmessage{This version of TeX cannot handle pictures.  Sorry.}}
\newhelp\disabledpictures
     {You will get a blank box in place of your picture.}
%
%
%
%
%
%
%
%
%
%
\def\FRONTPAGE{\ifvoid255\else\vfill\penalty-20000\fi
   \gl@bal\pagenumber=1     \gl@bal\chapternumber=0
   \gl@bal\equanumber=0     \gl@bal\sectionnumber=0
   \gl@bal\referencecount=0 \gl@bal\figurecount=0
   \gl@bal\tablecount=0     \gl@bal\frontpagetrue
   \gl@bal\lastf@@t=0       \gl@bal\footsymbolcount=0
   \gl@bal\cn@@false }

\def\papers{\papersize\headline=\paperheadline\footline=\paperfootline}
\def\papersize{\hsize=35pc \vsize=50pc \hoffset=0pc \voffset=1pc
   \advance\hoffset by\HOFFSET \advance\voffset by\VOFFSET
   \pagebottomfiller=0pc
   \skip\footins=\bigskipamount \normalspace }
\papers  
%
%
\newskip\lettertopskip       \lettertopskip=20pt plus 50pt
\newskip\letterbottomskip    \letterbottomskip=\z@ plus 100pt
\newskip\signatureskip       \signatureskip=40pt plus 3pt
\def\lettersize{\hsize=6.5in \vsize=8.5in \hoffset=0in \voffset=0.5in
   \advance\hoffset by\HOFFSET \advance\voffset by\VOFFSET
   \pagebottomfiller=\letterbottomskip
   \skip\footins=\smallskipamount \multiply\skip\footins by 3
   \singlespace }
\def\MEMO{\lettersize \headline=\letterheadline \footline={\hfil }%
   \let\rule=\memorule \FRONTPAGE \memohead }

\def\memodate{\afterassignment\MEMO \date }
\def\memit@m#1{\smallskip \hangafter=0 \hangindent=1in
    \Textindent{\caps #1}}
\def\subject{\memit@m{Subject:}}
\def\topic{\memit@m{Topic:}}
\def\from{\memit@m{From:}}
\def\to{\rel@x \ifmmode \rightarrow \else \memit@m{To:}\fi }
\def\memorule{\medskip\hrule height 1pt\bigskip}  
\def\memohead{\centerline{\fourteenrm MEMORANDUM}}
\newwrite\labelswrite
\newtoks\rw@toks
\def\letters{\lettersize
   \headline=\letterheadline \footline=\letterfootline
   \immediate\openout\labelswrite=\jobname.lab}

\let\letterhead=\rel@x
\def\addressee#1{\medskip\line{\hskip 0.75\hsize plus\z@ minus 0.25\hsize
                               \the\date \hfil }%
   \vskip \lettertopskip
   \ialign to\hsize{\strut ##\hfil\tabskip 0pt plus \hsize \crcr #1\crcr}
   \writelabel{#1}\medskip \noindent\hskip -\spaceskip \ignorespaces }
\def\rwl@begin#1\cr{\rw@toks={#1\crcr}\rel@x
   \immediate\write\labelswrite{\the\rw@toks}\futurelet\n@xt\rwl@next}
\def\rwl@next{\ifx\n@xt\rwl@end \let\n@xt=\rel@x
      \else \let\n@xt=\rwl@begin \fi \n@xt}
\let\rwl@end=\rel@x
\def\writelabel#1{\immediate\write\labelswrite{\noexpand\labelbegin}
     \rwl@begin #1\cr\rwl@end
     \immediate\write\labelswrite{\noexpand\labelend}}
\newtoks\FromAddress         \FromAddress={}
\newtoks\sendername          \sendername={}
\newbox\FromLabelBox
\newdimen\labelwidth          \labelwidth=6in
\def\makelabels{\afterassignment\Makelabels \sendername=}
\def\Makelabels{\FRONTPAGE \letterinfo={\hfil } \MakeFromBox
     \immediate\closeout\labelswrite  \input \jobname.lab\vfil\eject}
\let\labelend=\rel@x
\def\labelbegin#1\labelend{\setbox0=\vbox{\ialign{##\hfil\cr #1\crcr}}
     \MakeALabel }
\def\MakeFromBox{\gl@bal\setbox\FromLabelBox=\vbox{\Tenpoint
     \ialign{##\hfil\cr \the\sendername \the\FromAddress \crcr }}}
\def\MakeALabel{\vskip 1pt \hbox{\vrule \vbox{
        \hsize=\labelwidth \hrule\bigskip
        \leftline{\hskip 1\parindent \copy\FromLabelBox}\bigskip
        \centerline{\hfil \box0 } \bigskip \hrule
        }\vrule } \vskip 1pt plus 1fil }
\def\signed#1{\par \nobreak \bigskip \dt@pfalse \begingroup
  \everycr={\noalign{\nobreak
            \ifdt@p\vskip\signatureskip\gl@bal\dt@pfalse\fi }}%
  \tabskip=0.5\hsize plus \z@ minus 0.5\hsize
  \halign to\hsize {\strut ##\hfil\tabskip=\z@ plus 1fil minus \z@\crcr
          \noalign{\gl@bal\dt@ptrue}#1\crcr }%
  \endgroup \bigskip }
\newbox\letterb@x
\def\lettertext{\par \vskip\parskip \unvcopy\letterb@x \par }
\def\multiletter{\setbox\letterb@x=\vbox\bgroup
      \everypar{\vrule height 1\baselineskip depth 0pt width 0pt }
      \singlespace \topskip=\baselineskip }
\def\letterend{\par\egroup}
%
%
%
\newskip\frontpageskip
\newtoks\Pubnum   
\newtoks\Pubtype  \let\pubtype=\Pubtype
\newif\ifp@bblock  \p@bblocktrue
\def\PH@SR@V{\doubl@true \baselineskip=24.1pt plus 0.2pt minus 0.1pt
             \parskip= 3pt plus 2pt minus 1pt }
\def\PHYSREV{\papers\PhysRevtrue\PH@SR@V}
\let\physrev=\PHYSREV
\def\titlepage{\FRONTPAGE\papers\ifPhysRev\PH@SR@V\fi
   \ifp@bblock\p@bblock \else\hrule height\z@ \rel@x \fi }
\def\nopubblock{\p@bblockfalse}
\def\endpage{\vfil\break}
\frontpageskip=12pt plus .5fil minus 2pt
\Pubtype={}
\Pubnum={}
\def\p@bblock{\begingroup \tabskip=\hsize minus \hsize
   \baselineskip=1.5\ht\strutbox \topspace-2\baselineskip
   \halign to\hsize{\strut ##\hfil\tabskip=0pt\crcr
       \the\Pubnum\crcr\the\date\crcr\the\pubtype\crcr}\endgroup}
\def\title#1{\vskip\frontpageskip \titlestyle{#1} \vskip\headskip }
\def\author#1{\vskip\frontpageskip\titlestyle{\twelvecp #1}\nobreak}

\def\address#1{\par\kern 5pt\titlestyle{\twelvepoint\it #1}}
\def\andaddress{\par\kern 5pt \centerline{\sl and} \address}

\def\abstract{\par\dimen@=\prevdepth \hrule height\z@ \prevdepth=\dimen@
   \vskip\frontpageskip\centerline{\fourteenrm ABSTRACT}\vskip\headskip }

%
%
%
\def\ie{\hbox{\it i.e.}}       
\def\eg{\hbox{\it e.g.}}       
   
\def\\{\rel@x \ifmmode \backslash \else {\tt\char`\\}\fi }
\def\sequentialequations{\rel@x \if\equanumber<0 \else
  \gl@bal\equanumber=-\equanumber \gl@bal\advance\equanumber by -1 \fi }
\def\journal#1&#2(#3){\begingroup \let\journal=\dummyj@urnal
    \unskip, \sl #1\unskip~\bf\ignorespaces #2\rm
    (\afterassignment\j@ur \count255=#3), \endgroup\ignorespaces }
\def\j@ur{\ifnum\count255<100 \advance\count255 by 1900 \fi
          \number\count255 }
\def\dummyj@urnal{%
    \toks@={Reference foul up: nested \journal macros}%
    \errhelp={Your forgot & or ( ) after the last \journal}%
    \errmessage{\the\toks@ }}

\def\topspace{\hrule height 0pt depth 0pt \vskip}

\def\Buildrel#1\under#2{\mathrel{\mathop{#2}\limits_{#1}}}
\def\becomes#1{\mathchoice{\becomes@\scriptstyle{#1}}
   {\becomes@\scriptstyle{#1}} {\becomes@\scriptscriptstyle{#1}}
   {\becomes@\scriptscriptstyle{#1}}}
\def\becomes@#1#2{\mathrel{\setbox0=\hbox{$\m@th #1{\,#2\,}$}%
        \mathop{\hbox to \wd0 {\rightarrowfill}}\limits_{#2}}}

\let\int=\intop         
\def\lsim{\mathrel{\mathpalette\@versim<}}
\def\gsim{\mathrel{\mathpalette\@versim>}}
\def\@versim#1#2{\vcenter{\offinterlineskip
        \ialign{$\m@th#1\hfil##\hfil$\crcr#2\crcr\sim\crcr } }}
\def\big#1{{\hbox{$\left#1\vbox to 0.85\b@gheight{}\right.\n@space$}}}
\def\Big#1{{\hbox{$\left#1\vbox to 1.15\b@gheight{}\right.\n@space$}}}
\def\bigg#1{{\hbox{$\left#1\vbox to 1.45\b@gheight{}\right.\n@space$}}}
\def\Bigg#1{{\hbox{$\left#1\vbox to 1.75\b@gheight{}\right.\n@space$}}}
\def\){\mskip 2mu\nobreak }
%
%
%
\let\sec@nt=\sec
\def\sec{\rel@x\ifmmode\let\n@xt=\sec@nt\else\let\n@xt\section\fi\n@xt}
\def\obsolete#1{\message{Macro \string #1 is obsolete.}}
\def\firstsec#1{\obsolete\firstsec \section{#1}}
\def\firstsubsec#1{\obsolete\firstsubsec \subsection{#1}}
\def\thispage#1{\obsolete\thispage \gl@bal\pagenumber=#1\frontpagefalse}
\def\thischapter#1{\obsolete\thischapter \gl@bal\chapternumber=#1}
\def\splitout{\obsolete\splitout\rel@x}
\def\prop{\obsolete\prop \propto }
\def\nextequation#1{\obsolete\nextequation \gl@bal\equanumber=#1
   \ifnum\the\equanumber>0 \gl@bal\advance\equanumber by 1 \fi}
\def\BOXITEM{\afterassigment\B@XITEM\setbox0=}
\def\B@XITEM{\par\hangindent\wd0 \noindent\box0 }
%
%
%
%
%
%
   \def\unlock{\catcode`@=11}

   \def\lock{\catcode`@=12}

%
%
   \def\PRrefmark#1{\unskip~[#1]}
   \def\refitem#1{\ifPhysRev\r@fitem{[#1]}\else\r@fitem{#1.}\fi}
   \def\generatefootsymbol{%
      \ifcase\footsymbolcount\fd@f 13F \or \fd@f 279 \or \fd@f 27A
          \or \fd@f 278 \or \fd@f 27B
      \else%
         \ifnum\footsymbolcount <0%
            \xdef\footsymbol{\number-\footsymbolcount}%
         \else%
            \fd@f 203
               {\loop \ifnum\footsymbolcount >5
                  \fd@f{203 \footsymbol }
                  \advance\footsymbolcount by -1
                \repeat
               }
         \fi%
      \fi%
   }
   \def\OldPhysRevRefmark{\let\PRrefmark=\attach}
   \def\OldPRRefitem#1{\r@fitem{#1.}}
   \def\OldPhysRevRefitem{\let\refitem=\OldPRRefitem}
   \def\NPrefs{\let\refmark=\NPrefmark \let\refitem=\NPrefitem}
%
    \newif\iffileexists              \fileexistsfalse
    \newif\ifforwardrefson           \forwardrefsontrue
    \newif\ifamiga                   \amigafalse
    \newif\iflinkedinput             \linkedinputtrue
    \newif\iflinkopen                \linkopenfalse
    \newif\ifcsnameopen              \csnameopenfalse
    \newif\ifdummypictures           \dummypicturesfalse
    \newif\ifcontentson              \contentsonfalse
    \newif\ifcontentsopen            \contentsopenfalse
    \newif\ifmakename                \makenamefalse
    \newif\ifverbdone
    \newif\ifusechapterlabel         \usechapterlabelfalse
    \newif\ifstartofchapter          \startofchapterfalse
    \newif\iftableofplates           \tableofplatesfalse
    \newif\ifplatesopen              \platesopenfalse
    \newif\iftableoftables           \tableoftablesfalse
    \newif\iftableoftablesopen       \tableoftablesopenfalse
    \newif\ifwarncsname              \warncsnamefalse
%
    \newwrite\linkwrite
    \newwrite\csnamewrite
    \newwrite\contentswrite
    \newwrite\plateswrite
    \newwrite\tableoftableswrite
    \newread\testifexists
    \newread\verbinfile

    \newtoks\jobdir                  \jobdir={}
    \newtoks\tempnametoks            \tempnametoks={}
    \newtoks\oldheadline             \oldheadline={}
    \newtoks\oldfootline             \oldfootline={}
    \newtoks\subsectstyle            \subsectstyle={\Number}
    \newtoks\subsubsectstyle         \subsubsectstyle={\Number}
    \newtoks\runningheadlines        \runningheadlines={\relax}
    \newtoks\chapterformat           \chapterformat={\titlestyle}
    \newtoks\sectionformat           \sectionformat={\relax}
    \newtoks\subsectionformat        \subsectionformat={\relax}
    \newtoks\subsubsectionformat     \subsubsectionformat={\relax}
    \newtoks\chapterfontstyle        \chapterfontstyle={\bf}
    \newtoks\sectionfontstyle        \sectionfontstyle={\rm}
    \newtoks\subsectionfontstyle     \subsectionfontstyle={\rm}
    \newtoks\sectionfontstyleb       \sectionfontstyleb={\caps}
    \newtoks\subsubsectionfontstyle  \subsubsectionfontstyle={\rm}

    \newcount\subsectnumber           \subsectnumber=0
    \newcount\subsubsectnumber        \subsubsectnumber=0


   \newdimen\pictureindent           \pictureindent=15pt
   \newdimen\str
   \newdimen\squareht
   \newdimen\squarewd
   \newskip\doublecolskip
   \newskip\tableoftablesskip        \tableoftablesskip=\baselineskip


   \newbox\squarebox


   \newskip\sectionindent            \sectionindent=0pt
   \newskip\subsectionindent         \subsectionindent=0pt
  \def\thechapterhead{\relax}
  \def\thesectionhead{\relax}
  \def\thesubsecthead{\relax}
  \def\thesubsubsecthead{\relax}


   \def\GetIfExists #1 {
       \immediate\openin\testifexists=#1
       \ifeof\testifexists
           \immediate\closein\testifexists
       \else
         \immediate\closein\testifexists
         \input #1
       \fi
   }


   \def\stripbackslash#1#2*{\def\strippedname{#2}}

   \def\ifundefined#1{\expandafter\ifx\csname#1\endcsname\relax}

   \def\val#1{%
      \expandafter\stripbackslash\string#1*%
      \ifundefined{\strippedname}%
      \message{Warning! The control sequence \noexpand#1 is not defined.} ? %
      \else\csname\strippedname\endcsname\fi%
   }
%
%
   \def\CheckForOverWrite#1{%
      \expandafter\stripbackslash\string#1*%
      \ifundefined{\strippedname}%
      \else%
         \ifwarncsname
            \message{Warning! The control sequence \noexpand#1 is being
          overwritten.}%
          \else
          \fi
      \fi%
   }

   \def\FootNoteFonts{\Tenpoint}

   \def\Vfootnote#1{%
      \insert\footins%
      \bgroup%
         \interlinepenalty=\interfootnotelinepenalty%
         \floatingpenalty=20000%
         \singl@true\doubl@false%
         \FootNoteFonts%
         \splittopskip=\ht\strutbox%
         \boxmaxdepth=\dp\strutbox%
         \leftskip=\footindent%
         \rightskip=\z@skip%
         \parindent=0.5%
         \footindent%
         \parfillskip=0pt plus 1fil%
         \spaceskip=\z@skip%
         \xspaceskip=\z@skip%
         \footnotespecial%
         \Textindent{#1}%
         \footstrut%
         \futurelet\next\fo@t%
   }

   \def\csnamech@ck{%
       \ifcsnameopen%
       \else%
           \global\csnameopentrue%
           \immediate\openout\csnamewrite=\the\jobdir\jobname.csnames%
           \immediate\write\csnamewrite{\unlock}%
       \fi%
   }

   \def\linksch@ck{%
          \iflinkopen%
          \else%
              \global\linkopentrue%
              \immediate\openout\linkwrite=\the\jobdir\jobname.links%
          \fi%
   }

   \def\c@ntentscheck{%
       \ifcontentsopen%
       \else%
           \global\contentsopentrue%
           \immediate\openout\contentswrite=\the\jobdir\jobname.contents%
           \immediate\write\contentswrite{%
                \noexpand\titlestyle{Table of Contents}%
           }%
           \immediate\write\contentswrite{\noexpand\bigskip}%
       \fi%
   }

   \def\t@bleofplatescheck{%
       \ifplatesopen%
       \else%
           \global\platesopentrue%
           \immediate\openout\plateswrite=\the\jobdir\jobname.plates%
           \immediate\write\plateswrite{%
                \noexpand\titlestyle{Illustrations}%
           }%
           \immediate\write\plateswrite{%
              \unlock%
           }%
           \immediate\write\plateswrite{\noexpand\bigskip}%
       \fi%
   }

   \def\t@bleoftablescheck{%
       \iftableoftablesopen%
       \else%
           \global\tableoftablesopentrue%
          \immediate\openout\tableoftableswrite=\the\jobdir\jobname.tables%
           \immediate\write\tableoftableswrite{%
                \noexpand\titlestyle{Tables}%
           }%
           \immediate\write\tableoftableswrite{%
              \unlock%
           }%
           \immediate\write\tableoftableswrite{\noexpand\bigskip}%
       \fi%
   }


   \def\linkinput#1 {\input #1
       \iflinkedinput \relax \else \global\linkedinputtrue \fi
       \linksch@ck
       \immediate\write\linkwrite{#1}
   }


   \def\fil@#1 {%
       \ifdummypictures%
          \fileexistsfalse%
          \picfilename={}%
       \else%
          \fileexiststrue%
          \picfilename={#1}%
       \fi%
       \iflinkedinput%
          \iflinkopen \relax%
          \else%
            \global\linkopentrue%
            \immediate\openout\linkwrite=\the\jobdir\jobname.links%
          \fi%
          \immediate\write\linkwrite{#1}%
       \fi%
   }
   \def\Picture#1{%
      \gl@bal\advance\figurecount by 1%
      \CheckForOverWrite#1%
      \csnamech@ck%
      \immediate\write\csnamewrite{\def\noexpand#1{\the\figurecount}}%
      \xdef#1{\the\figurecount}\pl@t%
      \selfcaptionedtrue%
   }

   \def\s@vepicture{%
       \iffileexists \parsefilename \redopicturebox \fi%
       \ifdim\captionwidth>\z@ \else \captionwidth=\p@cwd \fi%
       \xdef\lastpicture{%
          \iffileexists%
             \setbox0=\hbox{\raise\the\yshift \vbox{%
                \moveright\the\xshift\hbox{\picturedefinition}}%
             }%
          \else%
             \setbox0=\hbox{}%
          \fi
          \ht0=\the\p@cht \wd0=\the\p@cwd \dp0=\the\p@cdp%
          \vbox{\hsize=\the\captionwidth \line{\hss\box0 \hss }%
          \ifcaptioned%
             \vskip\the\captionskip \noexpand\Tenpoint%
             \ifselfcaptioned%
                Figure~\the\figurecount.\enspace%
             \fi%
             \the\Caption%
          \fi }%
       }%
       \iftableofplates%
          \ifplatesopen%
          \else%
             \t@bleofplatescheck%
          \fi%
          \ifselfcaptioned%
             \immediate\write\plateswrite{%
                \noexpand\platetext{%
                \noexpand\item{\rm \the\figurecount .}%
                \the\Caption}{\the\pageno}%
             }%
          \else%
             \immediate\write\plateswrite{%
                \noexpand\platetext{\the\Caption}{\the\pageno}%
             }%
          \fi%
       \fi%
       \endgroup%
   }

   \def\platesout{%
      \ifplatesopen
         \immediate\closeout\plateswrite%
         \global\platesopenfalse%
      \fi%
      \input \jobname.plates%
      \lock%
   }

   \def\platetext#1#2{%
       \hbox to \hsize{\vbox{\hsize=.9\hsize #1}\hfill#2}%
       \vskip \tableoftablesskip \vskip\parskip%
   }


   \def\acksection#1{\par
      \ifnum\the\lastpenalty=30000\else \penalty-100\smallskip \fi
      \noindent\undertext{#1}\enspace \vadjust{\penalty5000}}

   \def\ack{\acksection{Acknowledgements:}}


   \def\pres@tpicture{%
       \gl@bal\linesabove=\linesabove
       \s@vepicture
       \setbox\picturebox=\vbox{
       \kern \linesabove\baselineskip \kern 0.3\baselineskip
       \lastpicture \kern 0.3\baselineskip }%
       \dimen@=\p@cht \dimen@i=\dimen@
       \advance\dimen@i by \pagetotal
       \par \ifdim\dimen@i>\pagegoal \vfil\break \fi
       \dimen@ii=\hsize
       \advance\dimen@ii by -\pictureindent \advance\dimen@ii by -\p@cwd
       \setbox0=\vbox to\z@{\kern-\baselineskip \unvbox\picturebox \vss }
   }

   \def\subspaces@t#1:#2;{%
      \baselineskip = \normalbaselineskip%
      \multiply\baselineskip by #1 \divide\baselineskip by #2%
      \lineskip = \normallineskip%
      \multiply\lineskip by #1 \divide\lineskip by #2%
      \lineskiplimit = \normallineskiplimit%
      \multiply\lineskiplimit by #1 \divide\lineskiplimit by #2%
      \parskip = \normalparskip%
      \multiply\parskip by #1 \divide\parskip by #2%
      \abovedisplayskip = \normaldisplayskip%
      \multiply\abovedisplayskip by #1 \divide\abovedisplayskip by #2%
      \belowdisplayskip = \abovedisplayskip%
      \abovedisplayshortskip = \normaldispshortskip%
      \multiply\abovedisplayshortskip by #1%
        \divide\abovedisplayshortskip by #2%
      \belowdisplayshortskip = \abovedisplayshortskip%
      \advance\belowdisplayshortskip by \belowdisplayskip%
      \divide\belowdisplayshortskip by 2%
      \smallskipamount = \skipregister%
      \multiply\smallskipamount by #1 \divide\smallskipamount by #2%
      \medskipamount = \smallskipamount \multiply\medskipamount by 2%
      \bigskipamount = \smallskipamount \multiply\bigskipamount by 4%
   }


   \def\makename#1{
       \global\makenametrue
       \global\tempnametoks={#1}
   }

   \def\nomakename#1{\relax}


   \def\savename#1{%
      \CheckForOverWrite{#1}%
      \csnamech@ck%
      \immediate\write\csnamewrite{\def\the\tempnametoks{#1}}%
   }

   \def\FootNoteFonts{\Tenpoint}

   \def\Vfootnote#1{%
      \insert\footins%
      \bgroup%
         \interlinepenalty=\interfootnotelinepenalty%
         \floatingpenalty=20000%
         \singl@true\doubl@false%
         \FootNoteFonts%
         \splittopskip=\ht\strutbox%
         \boxmaxdepth=\dp\strutbox%
         \leftskip=\footindent%
         \rightskip=\z@skip%
         \parindent=0.5%
         \footindent%
         \parfillskip=0pt plus 1fil%
         \spaceskip=\z@skip%
         \xspaceskip=\z@skip%
         \footnotespecial%
         \Textindent{#1}%
         \footstrut%
         \futurelet\next\fo@t%
   }
%

   \def\eqname#1{%
      \CheckForOverWrite{#1}%
      \rel@x{\pr@tect%
      \csnamech@ck%
      \ifnum\equanumber<0%
          \xdef#1{{\noexpand\f@m0(\number-\equanumber)}}%
          \immediate\write\csnamewrite{%
            \def\noexpand#1{\noexpand\f@m0 (\number-\equanumber)}}%
          \gl@bal\advance\equanumber by -1%
      \else%
          \gl@bal\advance\equanumber by 1%
          \ifusechapterlabel%
            \xdef#1{{\noexpand\f@m0(\ifcn@@ \chapterlabel.\fi%
               \number\equanumber)}%
            }%
          \else%
             \xdef#1{{\noexpand\f@m0(\ifcn@@%
                 {\the\chapterstyle{\the\chapternumber}}.\fi%
                 \number\equanumber)}}%
          \fi%
          \ifcn@@%
             \ifusechapterlabel
                \immediate\write\csnamewrite{\def\noexpand#1{(%
                  {\chapterlabel}.%
                  \number\equanumber)}%
                }%
             \else
                \immediate\write\csnamewrite{\def\noexpand#1{(%
                  {\the\chapterstyle{\the\chapternumber}}.%
                  \number\equanumber)}%
                }%
             \fi%
          \else%
              \immediate\write\csnamewrite{\def\noexpand#1{(%
                  \number\equanumber)}}%
          \fi%
      \fi}%
      #1%
   }

   \def\eqn{\eqno\eqname}

   \let\eqnalign=\eqname


   \def\APPENDIX#1#2{%
       \global\usechapterlabeltrue%
       \par\penalty-300\vskip\chapterskip%
       \spacecheck\chapterminspace%
       \chapterreset%
       \xdef\chapterlabel{#1}%
       \titlestyle{APPENDIX #2}%
       \nobreak\vskip\headskip \penalty 30000%
       \wlog{\string\Appendix~\chapterlabel}%
   }

   \def\REFNUM#1{%
      \CheckForOverWrite{#1} %
      \rel@x\gl@bal\advance\referencecount by 1%
      \xdef#1{\the\referencecount}%
      \csnamech@ck%
      \immediate\write\csnamewrite{\def\noexpand#1{\the\referencecount}}%
   }

   %

   \def\FIGNUM#1{
      \CheckForOverWrite{#1}%
      \rel@x\gl@bal\advance\figurecount by 1%
      \xdef#1{\the\figurecount}%
      \csnamech@ck%
      \immediate\write\csnamewrite{\def\noexpand#1{\the\figurecount}}%
   }


   \def\TABNUM#1{%
      \CheckForOverWrite{#1}%
      \rel@x \gl@bal\advance\tablecount by 1%
      \xdef#1{\the\tablecount}%
      \csnamech@ck%
      \immediate\write\csnamewrite{\def\noexpand#1{\the\tablecount}}%
   }


   \def\tableoftableson{%
      \global\tableoftablestrue%

      \gdef\TABLE##1##2{%
         \t@bleoftablescheck%
         \TABNUM ##1%
         \immediate\write\tableoftableswrite{%
            \noexpand\tableoftablestext{%
            \noexpand\item{\rm \the\tablecount .}%
                ##2}{\the\pageno}%
             }%
      }

      \gdef\Table##1{\TABLE\?{##1}Table~\?}
   }

   \def\tableoftablestext#1#2{%
       \hbox to \hsize{\vbox{\hsize=.9\hsize #1}\hfill#2}%
       \vskip \tableoftablesskip%
   }

   \def\tableoftablesout{%
      \iftableoftablesopen
         \immediate\closeout\tableoftableswrite%
         \global\tableoftablesopenfalse%
      \fi%
      \input \jobname.tables%
      \lock%
   }

%
%
%
%
%
%

   \def\contentsoff{\contentsonfalse}

   \def\f@m#1{\f@ntkey=#1\fam=\f@ntkey\the\textfont\f@ntkey\rel@x}
   \def\em@{\rel@x%
      \ifnum\f@ntkey=0\it%
      \else%
         \ifnum\f@ntkey=\bffam\it%
         \else\rm  %
         \fi%
      \fi%
   }

   \def\fontsoff{%
      \def\mit{\relax}%
      \let\oldstyle=\mit%
      \def\cal{\relax}%
      \def\it{\relax}%
      \def\sl{\relax}%
      \def\bf{\relax}%
      \def\tt{\relax}%
      \def\caps{\relax}%
      \let\cp=\caps%
   }


   \def\fontson{%
      \def\rm{\n@expand\f@m0}%
      \def\mit{\n@expand\f@m1}%
      \let\oldstyle=\mit%
      \def\cal{\n@expand\f@m2}%
      \def\it{\n@expand\f@m\itfam}%
      \def\sl{\n@expand\f@m\slfam}%
      \def\bf{\n@expand\f@m\bffam}%
      \def\tt{\n@expand\f@m\ttfam}%
      \def\caps{\n@expand\f@m\cpfam}%
      \let\cp=\caps%
   }

   \fontson
%


   \def\@alpha#1{\count255='140 \advance\count255 by #1\char\count255}
   \def\alphabetic{\@alpha}
   \def\@Alpha#1{\count255='100 \advance\count255 by #1\char\count255}
   \def\Alphabetic{\@Alpha}
   \def\@Roman#1{\uppercase\expandafter{\romannumeral #1}}
   \def\Roman{\@Roman}
   \def\@roman#1{\romannumeral #1}
   \def\roman{\@roman}
   \def\@number#1{\number #1}
   \def\Number{\@number}

   \def\leaderfill{\leaders\hbox to 1em{\hss.\hss}\hfill}

   \def\chapterinfo#1{%
      \line{%
         \ifcn@@%
            \hbox to \itemsize{\hfil\chapterlabel .\quad\ }%
         \fi%
         \noexpand{#1}\leaderfill\the\pagenumber%
      }%
   }

   \def\sectioninfo#1{%
      \line{%
         \ifcn@@%
            \hbox to 2\itemsize{\hfil\sectlabel \quad}%
          \else%
            \hbox to \itemsize{\hfil\quad}%
          \fi%
          \ \noexpand{#1}%
          \leaderfill \the\pagenumber%
      }%
   }

   \def\subsectioninfo#1{%
      \line{%
         \ifcn@@%
            \hbox to 3\itemsize{\hfil \quad\subsectlabel\quad}%
         \else%
            \hbox to 2\itemsize{\hfil\quad}%
         \fi%
          \ \noexpand{#1}%
          \leaderfill \the\pagenumber%
      }%
   }

   \def\subsubsecinfo#1{%
      \line{%
         \ifcn@@%
            \hbox to 4\itemsize{\hfil\subsubsectlabel\quad}%
         \else%
            \hbox to 3\itemsize{\hfil\quad}%
         \fi%
         \ \noexpand{#1}\leaderfill \the\pagenumber%
      }%
   }

   \def\CONTENTS#1;#2{
       {\let\makename=\nomakename
        \if#1C
            \immediate\write\contentswrite{\chapterinfo{#2}}%
        \else\if#1S
                \immediate\write\contentswrite{\sectioninfo{#2}}%
             \else\if#1s
                     \immediate\write\contentswrite{\subsectioninfo{#2}}%
                  \else\if#1x
                          \immediate\write\contentswrite{%
                              \subsubsecinfo{#2}}%
                       \fi
                  \fi
             \fi
        \fi
       }
   }

   \def\chapterreset{\gl@bal\advance\chapternumber by 1%
       \ifnum\equanumber<0 \else\gl@bal\equanumber=0 \fi%
       \gl@bal\sectionnumber=0  \gl@bal\let\sectlabel=\rel@x%
       \gl@bal\subsectnumber=0   \gl@bal\let\subsectlabel=\rel@x%
       \gl@bal\subsubsectnumber=0 \gl@bal\let\subsubsectlabel=\rel@x%
       \ifcn@%
           \gl@bal\cn@@true {\pr@tect\xdef\chapterlabel{%
           {\the\chapterstyle{\the\chapternumber}}}}%
       \else%
           \gl@bal\cn@@false \gdef\chapterlabel{\rel@x}%
       \fi%
       \gl@bal\startofchaptertrue%
   }

   \def\chapter#1{\par \penalty-300 \vskip\chapterskip%
       \spacecheck\chapterminspace%
       \gdef\thechapterhead{#1}%
       \gdef\thesectionhead{\relax}%
       \gdef\thesubsecthead{\relax}%
       \gdef\thesubsubsecthead{\relax}%
       \chapterreset \the\chapterformat{\the\chapterfontstyle%
          \ifcn@@\chapterlabel.~~\fi #1}%
       \nobreak\vskip\headskip \penalty 30000%
       {\pr@tect\wlog{\string\chapter\space \chapterlabel}}%
       \ifmakename%
           \csnamech@ck
           \ifcn@@%
              \immediate\write\csnamewrite{\def\the\tempnametoks{%
                 {\the\chapterstyle{\the\chapternumber}}}%
              }%
            \fi%
            \global\makenamefalse%
       \fi%
       \ifcontentson%
          \c@ntentscheck%
          \CONTENTS{C};{#1}%
       \fi%
       }%

   \def\section#1{\par \ifnum\lastpenalty=30000\else%
       \penalty-200\vskip\sectionskip \spacecheck\sectionminspace\fi%
       \gl@bal\advance\sectionnumber by 1%
       \gl@bal\subsectnumber=0%
       \gl@bal\let\subsectlabel=\rel@x%
       \gl@bal\subsubsectnumber=0%
       \gl@bal\let\subsubsectlabel=\rel@x%
       \gdef\thesectionhead{#1}%
       \gdef\thesubsecthead{\relax}%
       \gdef\thesubsubsecthead{\relax}%
       {\pr@tect\xdef\sectlabel{\ifcn@@%
          {\the\chapterstyle{\the\chapternumber}}.%
          {\the\sectionstyle{\the\sectionnumber}}\fi}%
       \wlog{\string\section\space \sectlabel}}%
       \the\sectionformat{\noindent\the\sectionfontstyle%
            {\ifcn@@\unskip\hskip\sectionindent\sectlabel~~\fi%
                \the\sectionfontstyleb#1}}%
       \par%
       \nobreak\vskip\headskip \penalty 30000%
       \ifmakename%
           \csnamech@ck%
           \ifcn@@%
              \immediate\write\csnamewrite{\def\the\tempnametoks{%
                 {\the\chapterstyle{\the\chapternumber}.%
                  \the\sectionstyle{\the\sectionnumber}}}
              }%
            \fi%
            \global\makenamefalse%
       \fi%
       \ifcontentson%
          \c@ntentscheck%
          \CONTENTS{S};{#1}%
       \fi%
   }

   \def\subsection#1{\par \ifnum\lastpenalty=30000\else%
       \penalty-200\vskip\sectionskip \spacecheck\sectionminspace\fi%
       \gl@bal\advance\subsectnumber by 1%
       \gl@bal\subsubsectnumber=0%
       \gl@bal\let\subsubsectlabel=\rel@x%
       \gdef\thesubsecthead{#1}%
       \gdef\thesubsubsecthead{\relax}%
       {\pr@tect\xdef\subsectlabel{\the\subsectionfontstyle%
           \ifcn@@{\the\chapterstyle{\the\chapternumber}}.%
           {\the\sectionstyle{\the\sectionnumber}}.%
           {\the\subsectstyle{\the\subsectnumber}}\fi}%
           \wlog{\string\section\space \subsectlabel}%
       }%
       \the\subsectionformat{\noindent\the\subsectionfontstyle%
         {\ifcn@@\unskip\hskip\subsectionindent%
          \subsectlabel~~\fi#1}}%
       \par%
       \nobreak\vskip\headskip \penalty 30000%
       \ifmakename%
           \csnamech@ck%
           \ifcn@@%
              \immediate\write\csnamewrite{\def\the\tempnametoks{%
                 {\the\chapterstyle{\the\chapternumber}}.%
                 {\the\sectionstyle{\the\sectionnumber}}.%
                 {\the\subsectstyle{\the\subsectnumber}}}%
              }%
            \fi%
            \global\makenamefalse%
       \fi%
       \ifcontentson%
          \c@ntentscheck%
          \CONTENTS{s};{#1}%
       \fi%
   }

   \def\subsubsection#1{\par \ifnum\lastpenalty=30000\else%
       \penalty-200\vskip\sectionskip \spacecheck\sectionminspace\fi%
       \gl@bal\advance\subsubsectnumber by 1%
       \gdef\thesubsubsecthead{#1}%
       {\pr@tect\xdef\subsubsectlabel{\the\subsubsectionfontstyle\ifcn@@%
           {\the\chapterstyle{\the\chapternumber}}.%
           {\the\sectionstyle{\the\sectionnumber}}.%
           {\the\subsectstyle{\the\subsectnumber}}.%
           {\the\subsubsectstyle{\the\subsubsectnumber}}\fi}%
           \wlog{\string\section\space \subsubsectlabel}%
       }%
       \the\subsubsectionformat{\the\subsubsectionfontstyle%
          \noindent{\ifcn@@\unskip\hskip\subsectionindent%
            \subsubsectlabel~~\fi#1}}%
       \par%
       \nobreak\vskip\headskip \penalty 30000%
       \ifmakename%
           \csnamech@ck%
           \ifcn@@%
              \immediate\write\csnamewrite{\def\the\tempnametoks{%
                {\the\chapterstyle{\the\chapternumber}.%
                 \the\sectionstyle{\the\sectionnumber}.%
                 \the\subsectionstyle{\the\subsectnumber}.%
                 \the\subsubsectstyle{\the\subsubsectnumber}}}%
              }%
            \fi%
            \global\makenamefalse%
       \fi%
       \ifcontentson%
          \c@ntentscheck%
          \CONTENTS{x};{#1}%
       \fi%
   }%

   \def\contentsinput{%
       \ifcontentson%
           \contentsopenfalse%
           \immediate\closeout\contentswrite%
           \global\oldheadline=\headline%
           \global\headline={\hfill}%
           \global\oldfootline=\footline%
           \global\footline={\hfill}%
           \fontsoff \unlock%
           \input \the\jobdir\jobname.contents%
           \fontson%
           \lock%
           \endpage%
           \global\headline=\oldheadline%
           \global\footline=\oldfootline%
       \else%
           \relax%
       \fi%
   }


       \def\phyzzxfootline{
           \footline={\ifletterstyle\the\letterfootline%
               \else\the\paperfootline\fi}%
       }

%

   {\obeyspaces}

   \def\verbfile#1{
       {\catcode`\\=12\catcode`\{=12
       \catcode`\}=12\catcode`\$=12\catcode`\&=12
       \catcode`\#=12\catcode`\%=12\catcode`\~=12
       \catcode`\_=12\catcode`\^=12\obeyspaces\obeylines\tt
       \verbdonetrue\openin\verbinfile=#1
       \loop\read\verbinfile to \inline
           \ifeof\verbinfile
               \verbdonefalse
           \else
              \leftline{\inline}
           \fi
       \ifverbdone\repeat
       \closein\verbinfile}
   }

   \def\boxit#1{\vbox{\hrule\hbox{\vrule\kern3pt%
       \vbox{\kern3pt#1\kern3pt}\kern3pt\vrule}\hrule}%
   }

   \def\square{%
      \setbox\squarebox=\boxit{\hbox{\phantom{x}}}
      \squareht = 1\ht\squarebox
      \squarewd = 1\wd\squarebox
      \vbox to 0pt{
          \offinterlineskip \kern -.9\squareht
          \hbox{\copy\squarebox \vrule width .2\squarewd height .8\squareht
              depth 0pt \hfill
          }
          \hbox{\kern .2\squarewd\vbox{%
            \hrule height .2\squarewd width \squarewd}
          }
          \vss
      }
   }

   \def\fboxit#1#2{
       \vbox{\hrule height #1
           \hbox{\vrule width #1
               \kern3pt \vbox{\kern3pt#2\kern3pt}\kern3pt \vrule width #1
           }
           \hrule height #1
       }
   }

   \let\eqnameold=\eqname

   \def\draft{\def\eqname##1{\eqnameold##1:{\tt\string##1}}
      \let\eqnalign = \eqname
   }
%
%
   \def\runningrightheadline{%
       \hfill%
       \tenit%
       \ifstartofchapter%
          \global\startofchapterfalse%
       \else%
          \ifcn@@ \the\chapternumber.\the\sectionnumber\quad\fi%
              {\fontsoff\thesectionhead}%
       \fi%
       \qquad\twelverm\folio%
   }

   \def\runningleftheadline{%
      \twelverm\folio\qquad%
      \tenit%
      \ifstartofchapter%
          \global\startofchapterfalse%
      \else%
         \ifcn@@%
             Chapter \the\chapternumber \quad%
         \fi%
         {\fontsoff\thechapterhead}%
         \hfill%
      \fi%
   }

   \runningheadlines={%
      \ifodd\pageno%
         \runningrightheadline%
      \else%
         \runningleftheadline%
      \fi
   }

%
%
%
%
%

   \font\dfont=cmr10 scaled \magstep5


   \newbox\cstrutbox
   \newbox\dlbox
   \newbox\vsk

   \setbox\cstrutbox=\hbox{\vrule height10.5pt depth3.5pt width\z@}

   \def\cstrut{\relax\ifmmode\copy\cstrutbox\else\unhcopy\cstrutbox\fi}

   \def\dl #1{\noindent\strut
       \setbox\dlbox=\hbox{\dfont #1\kern 2pt}%
       \setbox\vsk=\hbox{(}%
       \hangindent=1.1\wd\dlbox
       \hangafter=-2
       \strut\hbox to 0pt{\hss\vbox to 0pt{%
         \vskip-.75\ht\vsk\box\dlbox\vss}}%
       \noindent
   }

%
%

   \newdimen\fullhsize

   \fullhsize=6.5in
   \def\fullline{\hbox to\fullhsize}
   \let\l@r=L

   \newbox\leftcolumn
   \newbox\midcolumn

   \def\twocols{\hsize = 3.1in%
%
%
%
%
%
      \doublecolskip=.3333em plus .3333em minus .1em
      \global\spaceskip=\doublecolskip%
      \global\hyphenpenalty=0
      \singlespace
      \gdef\makeheadline{%
          \vbox to 0pt{ \skip@=\topskip%
          \advance\skip@ by -12pt \advance\skip@ by -2\normalbaselineskip%
          \vskip\skip@%
          \fullline{\vbox to 12pt{}\the\headline}\vss}\nointerlineskip%
      }%
      \def\makefootline{\baselineskip = 1.5\normalbaselineskip
           \fullline{\the\footline}
      }
      \output={%
          \if L\l@r%
             \global\setbox\leftcolumn=\columnbox \global\let\l@r=R%
          \else%
              \doubleformat \global\let\l@r=L%
          \fi%
          \ifnum\outputpenalty>-20000 \else\dosupereject\fi%
      }
      \def\doubleformat{
          \shipout\vbox{%
             \makeheadline%
             \fullline{\box\leftcolumn\hfil\columnbox}%
             \makefootline%
          }%
          \advancepageno%
      }
      \def\columnbox{\leftline{\pagebody}}
      \outer\def\twobye{%
          \par\vfill\supereject\if R\l@r \null\vfill\eject\fi\end%
      }%
   }

   \def\threecols{
       \hsize = 2.0in \tenpoint

      \doublecolskip=.3333em plus .3333em minus .1em
      \global\spaceskip=\doublecolskip%
      \global\hyphenpenalty=0

       \singlespace

       \def\makeheadline{\vbox to 0pt{ \skip@=\topskip
           \advance\skip@ by -12pt \advance\skip@ by -2\normalbaselineskip
           \vskip\skip@ \fullline{\vbox to 12pt{}\the\headline} \vss
           }\nointerlineskip
       }
       \def\makefootline{\baselineskip = 1.5\normalbaselineskip
                 \fullline{\the\footline}
       }

       \output={
          \if L\l@r
             \global\setbox\leftcolumn=\columnbox \global\let\l@r=M
          \else \if M\l@r
                   \global\setbox\midcolumn=\columnbox
                   \global\let\l@r=R
                \else \tripleformat \global\let\l@r=L
                \fi
          \fi
          \ifnum\outputpenalty>-20000 \else\dosupereject\fi
       }

       \def\tripleformat{
           \shipout\vbox{
               \makeheadline
               \fullline{\box\leftcolumn\hfil\box\midcolumn\hfil\columnbox}
               \makefootline
           }
           \advancepageno
       }

       \def\columnbox{\leftline{\pagebody}}

       \outer\def\threebye{
           \par\vfill\supereject
           \if R\l@r \null\vfill\eject\fi
           \end
       }
   }


%
%
%


   \everyjob{%
      \xdef\today{\monthname~\number\day, \number\year}
      
       \immediate\openin\testifexists=m
       \ifeof\testifexists
           \immediate\closein\testifexists
       \else
         \immediate\closein\testifexists
         \input m
       \fi
   yphyx.tex
      \ifforwardrefson%
         
       \immediate\openin\testifexists=\the
       \ifeof\testifexists
           \immediate\closein\testifexists
       \else
         \immediate\closein\testifexists
         \input \the
       \fi
   \jobdir\jobname.csnames
      \fi%
   }

\contentsoff

%
%
\def\phyzzx{PHY\setbox0=\hbox{Z}\copy0 \kern-0.5\wd0 \box0 X}
        
\message{ by V.K. and M.W. }

       \immediate\openin\testifexists=p
       \ifeof\testifexists
           \immediate\closein\testifexists
       \else
         \immediate\closein\testifexists
         \input p
       \fi
   hyzzx.local
\lock
%
%
%
\expandafter\def\expandafter\XXXX
	\expandafter{\csname\jobname\endcsname}
\def\YYYY{\phyzzx}
\ifx\XXXX\YYYY \let\next=\dump
\else \let\next=\relax \the\everyjob
\fi \next

\catcode `\@=11
\newskip\frontpageskip
\newtoks\Ftuvnum   \let\ftuvnum=\Ftuvnum	
\newtoks\Ificnum   \let\ificnum=\Ificnum
\newtoks\Pubtype  \let\pubtype=\Pubtype
\newif\ifp@bblock  \p@bblocktrue
\def\PH@SR@V{\doubl@true \baselineskip=15.0pt plus 0.1pt minus 0.1pt
             \parskip= 3pt plus 2pt minus 1pt \hsize=6in \vsize=8.6in}
\def\PHYSREV{\papers\PhysRevtrue\PH@SR@V}
\let\physrev=\PHYSREV
\def\titlepage{\FRONTPAGE\papers\ifPhysRev\PH@SR@V\fi
   \ifp@bblock\p@bblock \else\hrule height\z@ \rel@x \fi }
\def\nopubblock{\p@bblockfalse}
\def\endpage{\vfil\break}
\def\chapter#1{\par \penalty-300 \vskip\chapterskip
   \spacecheck\chapterminspace
   \chapterreset {\twelvebf{\noindent\ifcn@@\chapterlabel.~\fi #1\hfill}}
   \nobreak\vskip\headskip \penalty 30000
   {\pr@tect\wlog{\string\chapter\space \chapterlabel}} }
\frontpageskip=12pt plus .5fil minus 2pt
\Pubtype={}
\Ftuvnum={}
\Ificnum={}
\def\p@bblock{\begingroup \tabskip=\hsize minus \hsize
   \baselineskip=1.5\ht\strutbox \topspace-2\baselineskip
   \halign to\hsize{\strut ##\hfil\tabskip=0pt\crcr
       \the\Ftuvnum\crcr\the\Ificnum\crcr
	\the\date\crcr\the\pubtype\crcr}\endgroup}

\newcount\anni
\anni=\year
\advance\anni by -1900
\def\aapub{\afterassignment\aap@b\toks@}
\def\aap@b{\edef\n@xt{\Ftuvnum={FTUV\ \the\anni--\the\toks@}}\n@xt}

\def\bbpub{\afterassignment\bbp@b\toks@}
\def\bbp@b{\edef\n@xt{\Ificnum={IFIC\ \the\anni--\the\toks@}}\n@xt}

\let\ftuvnum=\aapub
\let\ificnum=\bbpub
\def\ftuvbin{
      \expandafter\ifx\csname binno\endcsname\relax%
         \expandafter\ifx\csname MailStop\endcsname\relax%
         \else%
            , Mail Stop \MailStop%
         \fi%
      \else%
         , Mail Stop \binno%
      \fi%
   }
\expandafter\ifx\csname eightrm\endcsname\relax
      \fi
\def\FTUV{\address{\tenit 
        Departamento de F\'{\i}sica Te\'{o}rica and IFIC,\break
        Centro Mixto Univ. de Valencia-CSIC\break 
	46100-Burjassot (Valencia), Spain\break
        {\tenrm azcarrag@evalvx.ific.uv.es, pbueno@lie.ific.uv.es}}}
\def\NAME{\centerline
{J.A. de Azc\'{a}rraga and J.C. P\'{e}rez Bueno}}

\def\pri{^{\, \prime }}
\VOFFSET=33pt
\papersize

\font\black=msbm10 scaled\magstep1

\def\bumpfootsymbolcount{\rel@x
   \iffrontpage \bumpfootsymbolpos \else \advance\lastf@@t by 1
     \ifPhysRev \bumpfootsymbolpos \else \bumpfootsymbolpos \fi \fi
   \gl@bal\lastf@@t=\pagen@ }
\def\title#1{\titlestyle{\twelvebf #1}}
\def\abstract{\par\dimen@=\prevdepth \hrule height\z@ \prevdepth=\dimen@
   \vskip\frontpageskip\centerline{\twelverm ABSTRACT}\vskip\headskip }
\def\ack{\par
      \ifnum\the\lastpenalty=30000\else \penalty-100\smallskip \fi
      \noindent{\bf Acknowledgements}:\quad\enspace \vadjust{\penalty5000}}

\def \apix {({B}.10)}
\def \apx {({B}.11)}
\def \apxi {({B}.12)}
\def \apxii {({B}.13)}
\def \apxiii {({B}.14)}
\def \apxiv {({B}.15)}
\def \apxv {({B}.16)}
\def \apxvi {({B}.17)}
\def \apxvii {({B}.18)}
\def \apxx {({B}.19)}

\def \apxxiii {({B}.22)}
\def \apxxiv {({B}.23)}

\def \apxxvi {({B}.25)}
\def \apxxvii {({B}.26)}
\def \apxxviii {({B}.27)}
\def \apxxviiia {({B}.28)}
\def \apxxix {({B}.29)}
\def \apxxx {({B}.30)}
\def \rotac {({\Number {2}}.1)}
\def \eucl {({\Number {2}}.2)}
\def \euclidean {({\Number {2}}.3)}
\def \galilei {({\Number {2}}.4)}
\def \VIi {({\Number {2}}.5)}
\def \VIii {({\Number {2}}.6)}
\def \redef {({\Number {3}}.1)}
\def \euclid {({\Number {3}}.2)}
\def \defbic {({\Number {3}}.3)}
\def \newgalbas {({\Number {3}}.4)}
\def \bicgalilei {({\Number {3}}.5)}
\def \defbicgalilei {({\Number {3}}.6)}
\def \VIiii {({\Number {3}}.7)}
\def \VIiv {({\Number {3}}.8)}
\def \VIv {({\Number {3}}.9)}
\def \VIvi {({\Number {3}}.10)}
\def \VIvii {({\Number {3}}.11)}
\def \VIvii {({\Number {3}}.12)}
\def \defua {({\Number {4}}.1)}
\def \rotacex {({\Number {4}}.2)}
\def \newrescal {({\Number {4}}.3)}
\def \heis {({\Number {4}}.4)}
\def \change {({\Number {4}}.5)}
\def \hei {({\Number {4}}.6)}
\def \acc {({\Number {4}}.7)}
\def \defcoci {({\Number {4}}.8)}
\def \newrotex {({\Number {4}}.9)}
\def \newres {({\Number {4}}.10)}
\def \rotexnew {({\Number {4}}.11)}
\def \changeii {({\Number {4}}.12)}
\def \rotexnewii {({\Number {4}}.13)}
\def \newbicross {({\Number {4}}.14)}
\def \heisa {({\Number {4}}.15)}
\def \heisab {({\Number {4}}.16)}
\def \insertx {({\Number {5}}.1)}
\def \insertxi {({\Number {5}}.2)}
\def \insertxii {({\Number {5}}.3)}
\def \insertxiii {({\Number {5}}.4)}
\def \inserti {({\Number {5}}.5)}
\def \insertii {({\Number {5}}.6)}
\def \hwrho {({\Number {5}}.7)}
\def \duali {({\Number {5}}.8)}
\def \insertiii {({\Number {5}}.9)}
\def \dualii {({\Number {5}}.10)}
\def \dualx {({\Number {5}}.11)}
\def \dualxi {({\Number {5}}.12)}
\def \plano {({\Number {6}}.1)}
\def \covariant {({\Number {6}}.2)}
\def \solut {({\Number {6}}.3)}
\def \defchi {({\Number {6}}.4)}
\def \solution {({\Number {6}}.5)}
\def \galplane {({\Number {6}}.6)}
\def \covargal {({\Number {6}}.7)}
\def \syst {({\Number {6}}.8)}
\def \galsolution {({\Number {6}}.9)}
\def \heiac {({\Number {6}}.10)}
\def \ecui {({\Number {6}}.11)}
\def \ecudos {({\Number {6}}.12)}
\def \ecutres {({\Number {6}}.13)}
\def \ecucuatro {({\Number {6}}.14)}
\def \difo {({A}.1)}
\def \difi {({A}.2)}
\def \difii {({A}.3)}
\def \difoo {({A}.4)}
\def \aapi {({B'}.1)}
\def \aapii {({B'}.2)}
\def \aapiiia {({B'}.3)}
\def \aapiii {({B'}.4)}
\def \aapiv {({B'}.5)}
\def \aapv {({B'}.6)}
\def \aapvi {({B'}.7)}
\def \aapvii {({B'}.8)}
\def \aapviii {({B'}.9)}
\def \aapix {({B'}.10)}
\def \aapx {({B'}.11)}
\def \aapxi {({B'}.12)}
\def \aapxii {({B'}.13)}
\def \aapxiii {({B'}.14)}
\def \aapxiv {({B'}.15)}
\def \aapxv {({B'}.16)}
\def \aapxvi {({B'}.17)}
\def \aapxvii {({B'}.18)}
\def \aapxx {({B'}.19)}
\def \aapxxi {({B'}.20)}
\def \aapxxii {({B'}.21)}
\def \aapxxiii {({B'}.22)}
\def \aapxxiv {({B'}.23)}
\def \aapxxv {({B'}.24)}
\def \aapxxvi {({B'}.25)}
\def \aapxxvii {({B'}.26)}
\def \aapxxviii {({B'}.27)}
\def \aapxxviiia {({B'}.28)}
\def \aapxxix {({B'}.29)}
\def \aapxxx {({B'}.30)}
\def \aapxxxi {({B'}.31)}
\def \aapxxxii {({B'}.32)}
\def \aapxxxiii {({B'}.33)}
\def \aapxxxx {({B'}.34)}
\def\JMP#1(#2){\journal J. Math. Phys. &#1(#2)}
\def\PHL#1(#2){\journal Phys. Lett. &#1(#2)}
\def\JA#1(#2){\journal J. Algebra &#1(#2)}
\def\LMP#1(#2){\journal Lett. Math. Phys. &#1(#2)}
\def\JPH#1(#2){\journal J. Phys. &#1(#2)}
\def\CMP#1(#2){\journal Commun. Math. Phys. &#1(#2)}

\REF\IW{E. \.In\"on\"u and E.P. Wigner, Proc. Nat. Acad. Sci 
{\bf 39} (1953), 510}

\REF\INONU{E. \.In\"on\"u, {\it Contractions of Lie groups and their 
representations} in {\it Group theor. concepts in elem. part. physics}, F. 
G\"ursey ed., Gordon and Beach, p. 391 (1964)}

\REF\AACON{V. Aldaya and J.A. de Azc\'{a}rraga, {\it Int. J. of Theor. Phys.}
{\bf 24} (1985), 141}

\REF\AI{J. A. de Azc\'{a}rraga and J. M. Izquierdo, 
{\it Lie algebras, Lie groups cohomology and some applications in physics}, 
Camb. Univ. Press (1995)}

\REF\SAL{E. J. Saletan \JMP 2 (61) 1}

\REF\DRI{V. G. Drinfel'd, in Proc. of the 1986 {\it Int. Congr. of Math.}, MSRI
Berkeley, vol {\bf I}, 798 (1987) (A. Gleason, ed.)}

\REF\JIM{M. Jimbo \LMP 10 (85) 63; ibid {\bf 11} (1986), 247}

\REF\FRT{L.D. Faddeev, N. Yu. Reshetikhin and L. A. Takhtajan, {\it Alg. i
Anal.} {\bf 1} (1989), 178 (Leningrad Math. J. {\bf 1} (1990), 193)}

\REF\FIRI{E. Celeghini, R. Giachetti, E. Sorace and M. Tarlini 
\JMP 31 (90) 2548; ibid {\bf 32} (1991), 1155, 1159}

\REF\CGSTa{E. Celeghini, R. Giachetti, E. Sorace, and M. Tarlini, {\it
Contractions of quantum groups}, in Lec. Notes Math. {\bf 1510}, (1992)
p. 221}

\REF\MB{S. Majid \JA 130 (90) 17;
{\it Israel J. Math.} {\bf 72} (1990), 133}

\REF\MAJSOB{S. Majid and Ya. S. Soibelman \JA 163 (94) 68}

\REF\SIN{W. Singer \JA 21 (72) 1}

\REF\BCM{R. J. Blattner, M. Cohen and S. Montgomery, {\it Trans. Am. Math.
Society} {\bf 298} (1986), 671; R. J. Blattner and S. Montgomery, {\it Pac. 
J. Math.} {\bf 137} (1989), 37}

\REF\MOL{R. Molnar \JA 47 (77) 29} 

\REF\LNRT {J. Lukierski, A. Nowicki, H. Ruegg and V.N. Tolstoy \PHL B264
(91) 331; 
J. Lukierski, H. Ruegg, and V.N. Tolstoy,
 {\it $\kappa$-quantum Poincar\'{e} 1994},
in {\it Quantum groups: formalism and applications}, J. Lukierski, Z. Popowicz
and J. Sobczyk eds, PWN (1994), p. 359}

\REF\MR{S. Majid and H. Ruegg \PHL B334 (94) 348} 

\REF\WOROa{S. L. Woronowicz \CMP 149 (92) 637}

\REF\SH{P. Schupp, P. Watts and B. Zumino \LMP 24 (92) 141}

\REF\BCGST{A. Ballesteros, E. Celeghini, R. Giachetti, E. Sorace and 
M. Tarlini \JPH A26 (93) 7495} 

\REF\BCHOS{A. Ballesteros, E. Celeghini, F. J. Herranz, M. A. del Olmo
and M. Santander \JPH A27 (94) L369}

\REF\GROMAN{N. A. Gromov and V. I. Man'ko \JMP 33 (92) 1374}

\REF\ElliSob{D. Ellinas and J. Sobczyk \JMP 36 (95) 1404}

\REF\BAM{W. K. Baskerville and S. Majid \JMP 34 (93) 3588}

\REF\BOGST{F. Bonechi, R. Giachetti, E. Sorace and M. Tarlini \CMP 169 (95) 627}

\REF\VK{L. L. Vaksman and L. I. Korogodskii, {\it Sov. Math. Dokl.} {\bf 39}
(1989), 173}

\REF\APBM{J. A. de Azc\'arraga, M. del Olmo, J. C. P\'erez Bueno, in 
preparation}

\REF\SIT{A. Sitarz \PHL B349 (95) 42}

\REF\APB{J. A. de Azc\'arraga and J. C. P\'erez Bueno \JMP 36 (95) 6879}

\REF\GLODZ{S. Giller, P. Kosi\'nski, M. Majewski, P. Ma\'slanka and J. Kunz 
\PHL B286 (92) 57}

\REF\MAC{A. J. Macfarlane \JPH A22 (89) 4581}

\REF\LBIE{L. C. Biedenharn \JPH A22 (89) L873}

\REF\AC{M. Arik and D. D. Coon \JMP 17 (76) 524}

\REF\KUL{P. P. Kulish, {\it Theor. Math. Phys.} {\bf 86} (1991), 108}

\REF\OS{C. H. Oh and K. Singh \JPH A27 (94) 5907}

\REF\QV{C. Quesne and N. Vansteenkiste, q-alg/9510001}

\REF\PM{P. Ma\'slanska \JMP 35 (94) 76}

\REF\SOB{J. Sobczyk, {\it Czech. J. Phys.} {\bf 46} (1996), 265}

\REF\HLR{V. Hussin, A. Lauzon and G. Rideau \LMP 31 (94) 159}

\REF\WORO{S. L. Woronowicz \CMP 122 (89) 125}

\catcode`\@=11

\font\tenmsa=msam10
\font\sevenmsa=msam7
\font\fivemsa=msam5
\font\tenmsb=msbm10
\font\sevenmsb=msbm7
\font\fivemsb=msbm5
\newfam\msafam
\newfam\msbfam
\textfont\msafam=\tenmsa  \scriptfont\msafam=\sevenmsa
  \scriptscriptfont\msafam=\fivemsa
\textfont\msbfam=\tenmsb  \scriptfont\msbfam=\sevenmsb
  \scriptscriptfont\msbfam=\fivemsb

\def\hexnumber@#1{\ifnum#1<10 \number#1\else
 \ifnum#1=10 A\else\ifnum#1=11 B\else\ifnum#1=12 C\else
 \ifnum#1=13 D\else\ifnum#1=14 E\else\ifnum#1=15 F\fi\fi\fi\fi\fi\fi\fi}

\def\msa@{\hexnumber@\msafam}
\def\msb@{\hexnumber@\msbfam}
\mathchardef\blacktriangleright="3\msa@49
\mathchardef\blacktriangleleft="3\msa@4A
\catcode`\@=\active

\def\R{\hbox{{\black R}}}
\newlist\mypoint=\alphabetic&)&0.5\itemsize;

\def\omm{\hat{\omega}}

\def\bic{\triangleright\!\!\!\blacktriangleleft}
\def\acti{\triangleleft}
\def\act{\triangleright}

\def\su{${\cal U}_q(su(2))$}
\def\U1{${\cal U}(u(1))$}
\def\UE2{${\cal U}_q({\cal E}(2))$}
\def\UEO{${\cal U}_{\omega}({\cal E}(2))$}

\def\H{{{\cal H}}}
\def\A{{{\cal A}}}
\def\rhoa{{\rho}}
\def\rhob{{\hat\rho}}
\def\FoE2{{\rm Fun}_\omega(E(2))}
\def\FqE2{{\rm Fun}_q(E(2))}
\def\FraE2{{\rm Fun}_\rhoa(\tilde E(2))}
\def\FrbE2{{\rm Fun}_\rhob(\tilde E(2))}
\def\FraHW{{\rm Fun}_\rhoa({\rm HW})}
\def\FrbHW{{\rm Fun}_\rhob({\rm HW})}
\def\tUEO{${\cal U}_{\rhoa}(\tilde{\cal E}(2))$}
\def\UOTR{{\cal U}_\omega({\cal T}r(2))}
\def\tEtO{${\tilde E}_\rhoa$}
\def\gal{{\cal U}_{\tilde\omega}({\cal G}(1+1))}
\def\galp{G^2_{\tilde\omega}}

\def\trrr{{\cal U}_{\tilde\omega}({\cal T}r(2))}

\def\UhoE2{${\cal U}_\rhob(\tilde{\cal E}(2))$}
\def\RL{\triangleright\!\!\!\blacktriangleleft}
\def\LR{\blacktriangleright\!\!\!\triangleleft}
\def\LLL{\triangleleft}
\def\RRR{\bar\triangleright}
\def\RIMO{\triangleright\!\!\!<}
\def\LECO{>\!\!\!\blacktriangleleft}
\def\LEMO{>\!\!\!\triangleleft}
\def\RICO{\blacktriangleright\!\!\!<}
\def\ABC#1{%
\global\usechapterlabeltrue%
\chapterreset%
\xdef\chapterlabel{#1}
}
\def\AA{H}
\def\HH{A}
\def\KK{K}
\def\balpha{\bar{\alpha}}
\def\bbeta{\bar{\beta}}
\def\bxi{\bar{\xi}}
\def\bpsi{\bar{\psi}}
\def\aa{h}
\def\bb{g}
\def\hh{a}
\def\gg{b}
\def\ff{c}
\def\EUC{{\cal U}({\cal E}(2))}
\physrev
\bumpfootsymbolpos

\ftuvnum={2}
\ificnum={3}
\date={January, 1996}
\titlepage
\vskip 1cm
\title{CONTRACTIONS, HOPF ALGEBRA EXTENSIONS AND\break
COVARIANT DIFFERENTIAL CALCULUS
\foot{To Jurek Lukierski on his 60th anniversary.}
}
\vskip 10pt
\NAME
\vskip 10pt
\FTUV
\vskip 10pt
\centerline{{\caps{Abstract}}}
{\leftskip 0.5in
\rightskip 0.5in
\baselineskip=13pt
{\tenrm We re-examine all the contractions related with 
the \su\ deformed algebra and study
the consequences that the contraction process has for their structure.
We also show using \su$\times$\U1\ as an example that,  
as in the undeformed case, the contraction may
generate Hopf algebra cohomology. We shall show that most of the different 
Hopf algebra
deformations obtained have a bicrossproduct or a cocycle bicrossproduct
structure, for which we shall also give their dual `group' versions.
The bicovariant differential calculi 
on the deformed spaces associated with the contracted algebras and the 
requirements for their existence are examined as well.}
\par
}
\chapter{Introduction}
As is well known, the standard Wigner-\.In\"on\"u contraction \refmark{\IW}
of simple Lie algebras with respect to a subalgebra leads to algebras which are 
the semidirect product of the preserved subalgebra and the resulting Abelian 
complement. Other types of contractions involving powers of the contraction 
parameter, first discussed in \refmark{\INONU}, 
may lead to a central extension structure.
Due to the singular nature of the contraction process, (non-simple) 
groups/algebras which are the
direct product/sum of two groups/algebras may not retain this direct product 
structure after the contraction limit
if the contraction affects suitably the central trivial 
extension; one may refer to these groups as being {\it pseudoextended} 
\refmark{\AACON,\AI}
when
the extension is trivial but behaves non-trivially under the contraction. 
A well known example is the direct product 
$P\times U(1)$, $P$ being the Poincar\'e group, for which a suitable limit 
leads to the centrally extended Galilei group 
\refmark{\SAL,\IW,\AACON,\AI}.

One of the interests of non-commutative geometry is to provide a rationale for
possible deformations of the spacetime manifold, which becomes a 
non-commutative algebra. 
By extending standard Lie group arguments about
quotient spaces, it 
is natural to associate these spacetime deformations with the 
deformation of inhomogeneous groups, which are non-simple.
Since the standard deformation procedure \refmark{\DRI, \JIM, \FRT} applies 
to the simple algebra/group 
case, the contraction of deformed simple algebras suggests
itself as a possible way of obtaining deformed inhomogeneous algebras.
This process usually requires involving the deformation 
parameter $q$ into the contraction \refmark{\FIRI, \CGSTa}, 
and is rather complicated; in fact, 
the contraction of deformed algebras/groups is besieged by the appearance
of divergences (the contraction is not always possible or the $R$-matrix 
diverges), and a complete theory is still lacking.
Clearly, the difficulty lies in having a well defined contraction process in 
{\it both} the algebra and coalgebra sectors.
 
Contraction is not, however, the only way of finding deformations of
inhomogeneous groups. 
Much in the same way we may construct Lie groups out of two by solving
the corresponding group extension problem (which always has a solution for 
Abelian kernel, precisely the semidirect extension, see \eg\ 
\refmark{\AI}),
we may look for a similar direct 
construction for Hopf algebras without thinking of obtaining them
by contraction.
Such a construction already exists for certain cases, and leads to the 
bicrossproduct and cocycle bicrossproduct 
structure of Hopf algebras of Majid \refmark{\MB,\MAJSOB} (see also
\refmark{\SIN,\BCM,\MOL};
a summary of Majid's theory is given in Appendix B).
For instance, the $\kappa$-Poincar\'e algebra of Lukierski {\it et al.}
\refmark{\LNRT}, which
is obtained from ${\cal U}_q(so(2,3))$ 
by a contraction involving the deformation parameter $q$ written 
as $q=\exp{1/\kappa R}$, where $R$ is the (de Sitter radius) 
contraction parameter so that $[\kappa]=L^{-1}$, 
has been shown {\refmark{\MR}} to possess 
such a bicrossproduct structure.
In this paper we intend to re-examine in this new light 
the simplest contraction examples, including the earliest ones 
\refmark{\CGSTa, \FIRI} 
(several of them discussed from various points of view in 
\refmark{\WOROa,\SH,\BCGST,\BCHOS,\GROMAN,\ElliSob,\allowbreak 
\BAM,\BOGST,\VK}).
We shall 
also look at the notion of central extension pseudocohomology for
Hopf algebras, and find that the contraction process generates Hopf algebra
extension cohomology as it does for its undeformed Lie counterpart.
We shall discuss both the `algebra' and `group' aspects of the deformed Hopf 
algebras, and study whether
they lead to a bicrossproduct or cocycle bicrossproduct structure. 
In contrast, the problems associated with the contraction 
and the $R$-matrix behaviour will not be discussed here. 
In fact, the constructions presented in sec. 5 may be considered as a way of 
avoiding the search for an $R$-matrix.
It would be interesting to perform a more general analysis of the consequences 
of the contraction process for the structure of the resulting deformed Hopf 
algebras.
We hope to report on this elsewhere \refmark{\APBM}.

The analysis of the differential calculus on the `spaces' associated with
the inhomogeneous deformed groups is also of importance; this has
been recently made for $\kappa$-spacetime algebras in 
\refmark{\SIT, \APB}.
It was shown there that the demand of covariance for the differential calculus 
required to enlarge the spacetime algebra by an element related 
to a central extension of the Hopf algebra; this phenomenon will also appear 
here for certain cases (sec. 6).
\chapter{Contractions of \su}
The well known \su\ deformed Hopf algebra 
is defined by ($q=e^z$)
$$
\eqalign{&[J_3,J_1]=J_2\quad,\quad[J_3,J_2]=-J_1\quad,\quad[J_1,J_2]={1
\over 2}[2J_3]_q={\sinh(2zJ_3)\over 2\sinh(z)}\quad;\cr
&\Delta J_{1,2}=\exp(-zJ_3)\otimes J_{1,2}+J_{1,2}\otimes\exp(zJ_3)\quad,
\quad \Delta J_3=J_3\otimes 1+1\otimes J_3\quad;\cr
&S(J_{1,2})=-\exp(zJ_3)J_{1,2}\exp(-zJ_3)\quad,\quad S(J_3)=-J_3\quad;\quad
\epsilon(J_{1,2,3})=0\quad .\cr}
\eqn\rotac
$$
Let us consider the different contractions of \su.

{\bf (1)}
\UE2.
The standard contraction procedure with respect the Hopf subalgebra generated
by $J_3$, implying the redefinitions $J_1=\epsilon^{-1}P_1\;,\;J_2=\epsilon
^{-1}P_2\;,\;J_3=J$, leads 
($\!\!$\refmark{\VK}; see also \refmark{\SH}) 
to
$$
\eqalign{&[J,P_1]=P_2\quad,\quad[J,P_2]=-P_1\quad,\quad[P_1,P_2]=0
\quad;\cr
&\Delta P_{1,2}=\exp(-zJ)\otimes P_{1,2}+P_{1,2}\otimes\exp(zJ)\quad,
\quad \Delta J=J\otimes 1+1\otimes J\quad;\cr
&S(P_{1,2})=-\exp(zJ)P_{1,2}\exp(-zJ)\quad,\quad S(J)=-J\quad;\quad
\epsilon(J,P_{1,2})=0\quad .\cr}
\eqn\eucl
$$
This deformation of the Euclidean algebra is a Hopf algebra where the 
deformation only appears at the coalgebra level, and will be denoted \UE2.  

{\bf (2)}
\UEO.
A second contraction, leading to another deformation 
of the Euclidean algebra, may be performed.
This contraction 
\refmark{\FIRI,\CGSTa}
requires writing previously $q=\exp(\epsilon\omega/2)$ since it is not 
performed with respect to a Hopf subalgebra
\foot{
It may be worth mentioning that contracting with respect a Hopf
subalgebra (as in the case (1) above) is not a sufficient condition to
define a contraction without involving the deformation parameter in it.
The above is, in fact, a rather exceptional case. 
}: 
it is performed with respect $J_2$, which is a Hopf subalgebra {\it only} for
$q=1$.
The redefinitions
$J_1=\epsilon^{-1}P_2\;,\;J_2=J\;,\;J_3=\epsilon^{-1}P_1\;,\;z=\epsilon\omega
/2$ in \rotac\ lead to \refmark{\FIRI,\CGSTa} the \UEO\ Euclidean Hopf
algebra
\foot{
If we want to look at $P_1$, $P_2$ as deformed translation generators,
$[P_i]=L^{-1}$, it is sufficient to take $[\epsilon]=L^{-1}\;,\;[\omega]=L$.
}
$$\eqalign{&[P_1,P_2]=0\quad,\quad[J,P_1]=P_2\quad,\quad[J,P_2]
=-{\sinh(\omega P_1)\over \omega}\quad;\cr
&\Delta P_1=P_1\otimes 1+1\otimes P_1\quad,\quad
\Delta P_2=\exp(-\omega P_1/2)\otimes P_2+P_2\otimes\exp(\omega P_1/2)\quad,
\cr 
&
\Delta J=\exp(-\omega P_1/2)\otimes J +J\otimes \exp(\omega P_1/2)
\quad;\quad S(P_{1,2})=-P_{1,2}\quad,\cr
&S(J)=-\exp(\omega P_1/2)J\exp(-\omega P_1/2)=
-J+{\omega\over 2}P_2\quad;
\quad\epsilon(J,P_{1,2})=0\quad .
\cr}
\eqn\euclidean
$$
\medskip
Besides the above, we may consider 
two `non-standard' contractions (\ie\ involving higher powers of the 
contraction parameter $\epsilon$). 
They are obtained by 
extending to the deformed case the generalized contraction in
\refmark{\INONU}.

{\bf (3)}
$\gal$.
A third contraction leads to a deformation of the Galilei algebra (the 
$(1+1)$ version of the $(1+3)$ deformed Galilei algebra in
\refmark{\GLODZ}).
We make the redefinitions 
\foot{
The parameter $\sigma\,,\,[\sigma]=TL^{-1/2}$, is introduced to give 
standard dimensions to the generators of the Galilei algebra 
($[\epsilon]=L^{-1/2}\,,\,[\tilde\omega]=T$),
but it disappears after the contraction.
}
$
z={\epsilon\tilde\omega\over \sigma}\,,\,
J_1=\sigma^{-1}\epsilon^{-1}\tilde V\,,\,
J_2=-\epsilon^{-2}\tilde X\,,\,
J_3=\epsilon^{-1}\sigma\tilde X_t.
$
By taking the limit $\epsilon\to 0$, we get
$$
\eqalign{
&
[\tilde X_t,\tilde V]=-\tilde X\quad,\quad
[\tilde X_t,\tilde X]=0\quad,\quad
[\tilde X,\tilde V]=0\quad;
\cr
&
\Delta \tilde X_t=\tilde X_t\otimes 1+1\otimes \tilde X_t\quad,\quad
\Delta \tilde X =\exp(-\tilde\omega \tilde X_t)\otimes \tilde X
+\tilde X\otimes\exp(\tilde\omega \tilde X_t)\quad,
\cr 
&
\Delta \tilde V =\exp(-\tilde\omega \tilde X_t)\otimes \tilde V
+\tilde V\otimes\exp(\tilde\omega \tilde X_t)
\quad;\quad 
S(\tilde X_t)=-\tilde X_t\quad,
\cr
&
S(\tilde X)=-\tilde X\quad,\quad
S(\tilde V)=-\exp(\tilde\omega \tilde X_t)\tilde V 
\exp(-\tilde\omega \tilde X_t)=-\tilde V+\tilde\omega\tilde X\quad,
\cr
&\epsilon(\tilde X_t,\tilde X,\tilde V)=0\quad .
\cr}
\eqn\galilei
$$
We will denote this deformed Galilei algebra by $\gal$. 
In the $\tilde\omega\to 0$ limit, eq. \galilei\ gives the Hopf 
structure of the enveloping algebra $\gal$ of the Galilei Lie algebra. 

{\bf (4)}
${\cal U}_{\hat\omega}(HW)$.
Finally, there is another contraction of \su. 
It is obtained
by making in \rotac\ the redefinitions $J_1=\epsilon^{-1}\bar X_q\;,\;
J_2=\epsilon^{-1}\bar X_p\;,\;J_3=\epsilon^{-2}\Xi$ and $z=\omm\epsilon^2/2$\  
\foot
{If one wishes to have $q$ and $p$ with dimensions of length and momentum (and 
$[X_q]=L^{-1},\,[X_p]=({\rm momentum})^{-1}$)
it is sufficient to modify the redefinitions to read 
$J_1=\epsilon^{-1}X_q\;,\;
J_2=\epsilon^{-1}\lambda X_p\;,\;J_3=\epsilon^{-2}\lambda \Xi\;,\; 
z=\omm\epsilon^2/2\lambda$, with $[\epsilon]=L^{-1}\;,\; [\lambda]=
[{\rm momentum}]L^{-1}\;,\break 
[\omm]={\rm action}\;,\;[\Xi]={\rm action}^{-1}$;
$\lambda$ disappears in the final expressions \VIi.}.
The result is the $\omm$-deformed Heisenberg-Weyl 
${\cal U}_{\omm}(HW)$ 
Hopf algebra
$$
\eqalign{
&
[\Xi,\bar X_q]=0\quad,\quad 
[\Xi,\bar X_p]=0\quad ,\quad
[\bar X_q,\bar X_p]={\sinh(\omm \Xi)\over \omm}\quad;
\cr
&
\Delta\bar X_{q,p}=\exp(-\omm\Xi/2)\otimes\bar X_{q,p}+
\bar X_{q,p}\otimes\exp(\omm\Xi/2)\quad ,\quad
\Delta\Xi=\Xi\otimes 1+1\otimes\Xi\quad;
\cr
&
S(\bar X_{q,p})=-\bar X_{q,p}\quad,\quad S(\Xi)=-\Xi
\quad;\quad\epsilon(\bar X_{q,p},\Xi)=0\quad,
\cr}
\eqn\VIi
$$
(denoted Heisenberg quantum group $H(1)_q$
in \refmark{\FIRI,\CGSTa}). 
By making the change of basis $X_{q,p}=\exp(-\omm\Xi/2)\bar X_{q,p}$ the 
${\cal U}_{\omm}(HW)$  
algebra takes the form
$$
\eqalign{
[\Xi,X_{p,q}]=0\quad,\quad 
[X_q,X_p]
&
={1-\exp(-2\omm\Xi)\over 2\omm}\quad;
\cr
\Delta X_{q,p}=X_{q,p}\otimes 1+\exp(-\omm\Xi)\otimes X_{q,p}\quad,\quad
&
\Delta\Xi=1\otimes\Xi+\Xi\otimes 1\quad;
\cr
\quad S(X_{q,p})=-X_{q,p}\exp(\omm\Xi)\quad;\quad
&
\epsilon(X_{q,p},\Xi)=0\quad;
\cr}
\eqn\VIii
$$
in the undeformed limit $\omm\rightarrow 0$, the standard 
expressions for the Hopf structure of ${\cal U}(HW)$ are recovered.

\chapter{Structure of the \su\ contractions}
As mentioned, the bicrossproduct 
\refmark{\MB,\MAJSOB} 
of Hopf algebras (see Appendix B) 
may be used as an alternative construction of deformed Hopf
algebras when the undeformed ones are not simple. 
Non-simple algebras may arise from  contraction, 
a process which for ordinary Lie algebras
leads to a semidirect product algebra.
Thus, it is worth exploring whether the above deformed Hopf 
algebras are the (right-left) bicrossproduct $\H\bic \A$ of two Hopf algebras 
$\H$ and $\A$ or have a cocycle bicrossproduct structure.
The notation $\H\bic \A$, for instance, indicates  
that $\A$ is a right $\H$-module algebra for the
{\it right}
action $\alpha:\A\otimes \H\to \A\;,\;\alpha(a,h)\equiv a\acti h\;,$ and that
$\H$ is a left $\A$-comodule coalgebra for the {\it left} 
coaction $\beta:\H\to
\A\otimes \H$ ($\H$ is a {\it left quantum space}); 
$\alpha$ and $\beta$ must also
satisfy certain compatibility conditions \refmark{\MB,\MAJSOB}.

{\bf (1)}
Let us first consider \UE2, eqs. \eucl. At the algebra level it has a 
semidirect structure. 
However, if we take $\A$ as the undeformed Hopf algebra generated 
by $P_1\;,\;P_2$ and $\H$ as that generated by $J$ we see that with 
independence of $\beta$, we cannot reproduce $\Delta(P_{1,2})$ in \UE2; in 
fact, $P_1\,,\, P_2$ in \eucl\ do not generate a Hopf subalgebra of \UE2.
Thus, \UE2 has not a bicrossproduct structure.

{\bf (2)}
Let us now look at \UEO, eqs. \euclidean. The redefinitions 
$$P_x=P_1\quad,\quad P_y=\exp(-\omega P_1/2)P_2\quad,\quad 
J^{\pri}=\exp(-\omega P_1/2)J\quad,\eqn\redef$$
allow us to write \UEO\ in terms of $(P_x,P_y,J^{\pri})$ in the form
$$
\eqalign{
[P_x,P_y]=0\quad,\quad[J^{\pri},P_x]=P_y\quad,
&
\quad
[J^{\pri},P_y]=
-{1\over 2\omega}\left(1-\exp(-2\omega P_x)\right)-{\omega\over 2}P_y^2\quad;
\cr
\Delta P_x=P_x\otimes 1+1\otimes P_x\quad,
&
\quad
\Delta P_y=P_y\otimes 1+\exp(-\omega P_x)\otimes P_y\quad,
\cr
\Delta J^{\pri}=J^{\pri}\otimes 1+
&
\exp(-\omega P_x)\otimes J^{\pri}\quad;\quad \epsilon(J,P_{x,y})=0\quad;
\cr
S(P_x)=-P_x\quad,\quad S(P_y)=
&
-\exp(\omega P_x)P_y\quad,\quad 
S(J^{\pri})=-\exp(\omega P_x)J^{\pri}\quad.
\cr}
\eqn\euclid
$$
If we now take for $\A$ the commutative non-cocommutative Hopf translation 
subalgebra $\UOTR$ of $(P_x,P_y)$ 
contained in \euclid\ and $\H$
is the commutative and cocommutative algebra generated by $J'$, 
the bicrossproduct structure $\H\bic \A$ of \euclid\ is
exhibited if
$$
\alpha(P_{x,y},J^{\pri})\equiv P_{x,y}\acti J^{\pri}:=[P_{x,y},J^{\pri}]\quad,
\quad\beta(J^{\pri}):=\exp(-\omega P_x)\otimes J^{\pri}\quad,
\eqn\defbic
$$
since it may be seen that the compatibility axioms 
\foot{
The formulae (B.-) refer to the corresponding (B'.-) ones given in Appendix B; 
they may be found in the original papers \refmark{\MB}
or in the Appendix of \refmark{\APB}
(there with the same numbering).
}
\apix, \apx, \apxi, 
\apxii\ and \apxiii\ are satisfied and that \apxv, \apxvi, \apxvii, define
the coproducts, antipodes and counits in \euclid. This shows that
\UEO$={\cal U}(u(1))\bic\UOTR$.

{\bf (3)}
Consider now the deformed $(1+1)$ 
Galilei Hopf algebra $\gal$ of \galilei. 
It was found in \refmark{\APB} (for the $(1+3)$ case) that it 
is also endowed with a bicrossproduct structure. 
To show this, we make the redefinitions
$$
X_t=\tilde X_t\quad,\quad
X=\exp(-\tilde\omega\tilde X_t)\tilde X\quad,\quad
V=\exp(-\tilde\omega\tilde X_t)\tilde V\quad.
\eqn\newgalbas
$$
With them, the $\gal$ Hopf algebra takes the form
$$
\eqalign{ 
&
[X_t,V]=-X\quad,\quad
[X,V]= {\tilde\omega}X^2\quad,\quad
[X,X_t]=0\quad;
\cr &
\Delta X_t=X_t\otimes 1+1\otimes X_t\quad,\quad
\Delta X=X\otimes 1+\exp(-2\tilde\omega X_t)\otimes X\quad,
\cr &
\Delta V=V\otimes 1+\exp(-2\tilde\omega X_t)\otimes V\quad;\quad 
\epsilon(V,X,X_t)=0\quad;
\cr &
S(X_t)=-X_t\quad,\quad S(X)=-\exp(2\tilde\omega X_t)X\quad,\quad 
S(V)=-\exp(2\tilde\omega X_t)V\quad,
\cr}
\eqn\bicgalilei
$$
(which is eq. (6.1) in \refmark{\APB} for ${\cal G}_{\tilde\kappa}$ with 
$1/2\tilde\kappa=\tilde\omega$).
The bicrossproduct structure is summarized in the definitions of the action 
$\alpha$ and the coaction $\beta$ ($\A$ is the Abelian, non-cocommutative
Hopf subalgebra $\trrr$ generated by 
$X$ and $X_t$, and $\H$ is given by the commutative and cocommutative 
Hopf algebra generated by $V$)
$$
\eqalign{
\alpha(X,V)\equiv X\acti V:=[X,V]=\tilde\omega X^2\quad,
&\quad
\alpha(X_t,V)\equiv X_t\acti V:=[X_t,V]=-X\quad,
\cr
\beta(V):=\exp
&
(-2\tilde\omega X_t)\otimes V\quad.
\cr}
\eqn\defbicgalilei
$$
It may be shown that the bicrossproduct conditions
are verified and hence that
$\gal=$\U1$\bic\trrr$. 
\foot{
Note that the expressions \galilei, \bicgalilei\ and \defbicgalilei\ may be 
obtained from standard contraction of their analogous ones in the Euclidean case
\euclidean, \euclid\ and \defbic.
} 

{\bf (4)}
Finally, we now show that ${\cal U}_{\omm}(HW)$ [\VIii] has {\it both} 
a bicrossproduct and a cocycle bicrossproduct 
structure. This parallels the fact that the Heisenberg-Weyl ${\rm (HW)}$ 
Lie group,
$(q^\prime,p^\prime,\theta^\prime)(q,p,\theta)=(q^\prime + q,p^\prime +p,
\,\theta^\prime +\theta 
+q^\prime p)$, may be considered as the semidirect extension of $\R\;({\rm 
coordinate}\; q)$ by the 
invariant subgroup $\R^2\;({\it id.}\; p,\theta)$ (the action 
of $\R$ on $\R^2$ being 
given by $q:(p,\theta)\mapsto(p,qp+\theta)$), {\it or} 
as a central extension of $\R^2
\;({\rm coordinates}\;p,q)$ by $\R\;({\it id.}\;\theta)$ 
(the $\R$-valued two-cocycle being given in its `asymmetric' form
$\xi(p^\prime,q^\prime;p,q)=q^\prime p$).

4a)
The {\it bicrossproduct structure} follows taking for $\A$ the Abelian
$\omm$-deformed Hopf 
subalgebra generated by $X_q$ and $\Xi$ in \VIii, for $\H$ the 
undeformed algebra generated by $X_p$ $(\Delta X_p=X_p\otimes 1+1\otimes 
X_p,\;S(X_p)=-X_p,\;\epsilon(X_p)=0)$ 
and for $\alpha$ and $\beta$ 
$$
X_q\acti X_p=[X_q,X_p]\quad,\quad\Xi\acti X_p=0\quad;\quad 
\beta(X_p)=\exp(-\omm\Xi)\otimes X_p\quad.
\eqn\VIiii
$$
This induces the appropriate coproduct for $X_p$ (identified as $X_p\otimes 1$ 
in $\H\otimes \A$) and antipode (eq. \VIii) from \apxv\ and \apxvii,  
respectively; the commutators in \VIii\ follow from \apxiv.

4b)
The {\it cocycle bicrossproduct structure}  
is constructed from the {\it undeformed} Hopf algebras 
$\A={\cal U}(\Xi)\;({\rm generated\; by}\; \Xi)$ 
and $\H=({\cal T}r(2))\; ({\it id.}\; X_q,X_p)$. 
Since we wish to obtain a deformation of a 
central extension algebra the action $\alpha$ must be trivial. We take
$$
\alpha(\Xi,X_{q,p})=[\Xi,X_{q,p}]=0\quad,
\quad \beta(X_{q,p})=\exp(-\omm\Xi)\otimes 
X_{q,p}
\eqn\VIiv
$$
and $\xi:\H\otimes \H\rightarrow \A\; {\rm (antisymmetric\; cocycle)},
\;\psi:\H\rightarrow \A\otimes \A$ given by
$$
\xi(X_q,X_p)=-\xi(X_p,X_q)={1\over 2}
\exp(-\omm\Xi){\sinh(\omm\Xi)\over\omm}\quad,\quad
\psi(X_{q,p})=0
\eqn\VIv
$$
(\ie\ $\psi$ trivial, $\psi(h)=1\otimes 1\epsilon(h)$) 
plus \apxx\ and \apxxviiia. 
Thus, the deformed character of the resulting algebra (and hence $\omm$) enters
in this case through $\beta$ and $\xi$ only. Since 
$\alpha$ and $\psi$ are trivial and $\A$ and $\H$ have the cocommutative Hopf 
algebra structure associated with the Abelian 
enveloping algebras ${\cal U}(\Xi ),\;{\cal U}({\cal T}r(2))$, it is not 
difficult to check that the compatibility conditions \apxiii, \apxxvi, 
\apxxvii\  
and \apxxviii\ are fulfilled. Moreover, \apxxix\ and \apxxx\ reduce using 
\apxx\ to
$$
(h\otimes a)(g\otimes b)=hg\otimes ab+1\otimes\xi(h,g)ab\quad,\quad h,g\in 
{\cal H}
\eqn\VIvi
$$
$$
\Delta(h\otimes 1)=h\otimes 1\otimes 1\otimes 1+1\otimes h^{(1)}\otimes 
h^{(2)}\otimes 1\quad,\quad h,g\in{\cal H}\quad.
\eqn\VIvii
$$
Denoting the elements $(h\otimes 1)$ and $(g\otimes 1)$ in $\H\otimes \A$ by 
$\tilde h$ and $\tilde g$, 
eq. \VIvi\ leads to $[\tilde h,\tilde g]=2\xi(h,g)\;([h,g]=0\;
{\rm in}\; \H)$ so that the commutators in 
\VIii\ are recovered for $\xi$ given by \VIv. Similarly,
$$\Delta(X_{q,p})=X_{q,p}\otimes 1\otimes 1\otimes 
1+1\otimes\exp(-\omm\Xi)\otimes X_{q,p}\otimes 1\equiv X_{q,p}\otimes 
1+\exp(-\omm\Xi)\otimes X_{q,p}\eqn\VIvii$$
plus $\Delta(\Xi)=\Xi\otimes 1+1\otimes\Xi$.
With $\epsilon(X_{q,p})=0=\epsilon(\Xi)$, the Hopf algebra 
structure of
${\cal U}_{\omm}(HW)$
is obtained by adding the antipode as defined by \VIii.

\chapter{The cocycle extended Euclidean Hopf algebras \tUEO, \UhoE2}
Consider \UEO$\times$\U1. This Hopf algebra has a trivial 
central factor and, as such, it
might have been obtained from \su$\times$\U1 by contraction,  since
the redefinitions given in sec. 2 (2) (and \redef) do not affect 
the \U1 part $(\Xi)$.
However a generalization of the pseudocohomology 
mechanism \refmark{\AACON,\AI} 
mentioned in the introduction may also be used 
here to obtain non-trivial extensions of
Hopf algebras by contracting trivial products (see \refmark{\APB} for
the case of the deformed extended (1+3) Galilei Hopf algebra). 
We now find two deformations of the centrally extended
Euclidean algebra using this procedure.

{\bf (a)}
\tUEO.
Consider the \su$\times$\U1\ Hopf algebra generated by 
$(J_1,J_2,\allowbreak J_3,\Xi')$
given by eqs. \rotac\ plus the
\U1\ relations 
$$
\Delta\Xi'=\Xi'\otimes 1+1\otimes\Xi'\quad,\quad
S(\Xi')=-\Xi'\quad,\quad\epsilon(\Xi')=0\quad;\quad
[\Xi',{\rm all}]=0\quad.
\eqn\defua
$$

The redefinition $J_3=J_3^{\pri}+\Xi^{\pri}$ \refmark{\FIRI} 
leaves \rotac\ and
\defua\ unchanged but for
$$\eqalign{[J_1,J_2]=&{\sinh(2z(J_3^{\pri}+\Xi^{\pri}))\over 2\sinh(z)}
\quad;\cr
\Delta J_{1,2}=\exp(-z(J_3^{\pri}+\Xi^{\pri}))\otimes & J_{1,2}+
J_{1,2}\otimes\exp(z(J_3^{\pri}+\Xi^{\pri}))\quad,\cr
\Delta J_3^{\pri}=&J_3^{\pri}\otimes 1+1\otimes J_3^{\pri}\quad;\cr
S(J_{1,2})=-\exp(z(J_3^{\pri}+\Xi^{\pri}))J_{1,2}&
\exp(-z(J_3^{\pri}+\Xi^{\pri}))
\quad,\quad 
S(J_3^{\pri})=-J_3^{\pri}\quad.\cr}\eqn\rotacex$$
Because $[J_1, J_2]$ involves $\Xi^{\pri}$, we refer to \su$\times$\U1\ in the
form \rotacex\ as a {\it pseudoextension} (the trivial direct product 
structure is disguised beneath the election of the generators).

To obtain a non-trivial Hopf algebra extension from it, we now make 
a rescaling {\it involving} $\Xi^{\pri}$,
$$
J_1=\epsilon^{-1}X_1\quad,\quad  J_2=\epsilon^{-1}X_2\quad,\quad 
J_3^{\pri}=N\quad,\quad \Xi^{\pri}=\Xi/\epsilon^2\quad,
\eqn\newrescal
$$
redefine $z$ as $z=\rhoa\epsilon^2$ and take the limit 
$\epsilon\to 0$. The resulting Hopf algebra is given by
$$
\eqalign{
&
[N,X_1]=X_2\ ,\quad 
[N,X_2]=-X_1\ ,\quad  
[X_1,X_2]={\sinh 2\rhoa\Xi\over 2\rhoa}\ ,\quad  
[\Xi,{\rm all}]=0\quad;
\cr
&
\Delta X_i=\exp(-{\rhoa}\Xi)\otimes X_i+ X_i\otimes
\exp({\rhoa}\Xi)\;(i=1,2)\ ,\ 
\Delta N=N\otimes 1+1\otimes N\ ;
\cr
&
S(X_{1,2})=-X_{1,2}\quad ,\quad 
S(N)=-N\quad ,\quad 
S(\Xi)=-\Xi\quad ;\quad  \epsilon(X_{1,2},N,\Xi)=0\ .
\cr}
\eqn\heis
$$
This Hopf algebra will be denoted by \tUEO.

It is convenient to make in \heis\ the change
$$ 
Y_i=\exp(-\rhoa\Xi)X_i\quad.\eqn\change$$
This modifies only
$$
\eqalign{
[Y_1,Y_2]=
&
{1\over 4\rhoa}(1-\exp(-4\rhoa\Xi))\quad,
\cr 
\Delta Y_i=\exp(-2\rhoa\Xi)\otimes Y_i+Y_i\otimes 1\quad,
&
\quad S(Y_i)=-\exp(2\rhoa\Xi)Y_i\quad,
\cr}
\eqn\hei
$$
which reproduces the Heisenberg-Weyl ${\cal U}_{\omm}(HW)$ 
Hopf algebra of \VIii\ with $2\rho=\omm$.

a1)
\tUEO\ has the bicrossproduct 
structure $\H\bic\A$, in which $\A$ is the deformed Heisenberg-Weyl 
${\cal U}_\rhoa(HW)$ Hopf subalgebra in \tUEO\ generated 
by $(Y_1,Y_2,\Xi)$ with primitive coproduct for $\Xi$ and $\Delta(Y_i)$ given 
in \hei, and $\H$ is the commutative and cocommutative algebra generated by 
$N$. The right action $\acti$ of $N$ on $\A$ is then designed to reproduce the 
commutators in \tUEO
$$
\alpha(Y_1,N)=[Y_1,N]=-Y_2\quad,\quad
\alpha(Y_2,N)=[Y_2,N]=Y_1,\quad,\quad 
\Xi\acti N=0\quad,
\eqn\acc
$$
and the coaction $\beta$ is taken to be trivial, $\beta(N)=1\otimes N$, since 
the coproducts in both $\H$ and $\A$ are already those in \tUEO.

a2) 
The cocycle extension structure of \tUEO\ 
is achieved by taking $\A$ generated by $\Xi$ [eq. \defua] and
$\H$ as the undeformed Euclidean algebra $\EUC$. The action
$\alpha$ of $\H$ on $\A$ is trivial (we want $\Xi$ to be central), 
and so is the map $\psi$ (\apxxiii, \apxxiv); 
the antisymmetric cocycle $\xi$ and coaction $\beta$ are given
by (cf. \hei)
$$
\xi(Y_1,Y_2)={1\over 8\rhoa}(1-\exp(-4\rhoa\Xi))\;,\;
\beta(Y_i)=\exp(-2\rhoa\Xi)\otimes Y_i\;,\;\beta(N)=1\otimes N
\eqn\defcoci
$$
(the coaction on $N$ is trivial).
We may check that all relations \apxx-\apxxvii, \apxxviii
\foot{
Since $\psi$ and $\LLL$ are trivial, this formula reduces to
$
\Delta\xi(h\otimes g)=\xi(h_{(1)}\otimes g_{(1)})h_{(2)}^{\ (1)}g_{(2)}^{\ 
(1)}\otimes\xi(h_{(2)}^{\ (2)}\otimes g_{(2)}^{\ (2)})\ .
$
},
\apxxviiia\ 
are fulfilled and that
\apxxix-\apxxx\ then reproduce \hei;
thus, 
\tUEO\ has a cocycle bicrossproduct structure.

\medskip
{\bf (b)}
\UhoE2.
Consider again the algebra ${\cal U}_q(su(2))\times {\cal U}(u(1))$ given by 
eqs. \rotac\ plus the relations \defua\ for the central $u(1)$ generator, 
now denoted $\hat\Xi'$. 
The redefinition $J_1=J_1'+\hat\Xi'$ leaves \defua\ and 
\rotac\ unchanged but for
$$
\eqalign{
[J_3,J_2]=-(J_1'+\hat\Xi')&\quad,\cr
\Delta J_1'=
\exp({-zJ_3})\otimes J_1'+J_1'\otimes\exp({zJ_3})
&
\cr
+(\exp({-zJ_3})-1)\otimes\hat\Xi'+\hat\Xi'\otimes
&
(\exp({zJ_3})-1)\quad.
\cr}
\eqn\newrotex
$$
If we now make the rescaling
$$
J_3=\epsilon^{-1}P_x\quad,\quad 
J_2=-\epsilon^{-1}P_y\quad,\quad 
J_1'=J'\quad,\quad
\hat\Xi'=\epsilon^{-2}\hat\Xi\quad,
\eqn\newres
$$
and set $z=\epsilon^3\rhob$, in the limit $\epsilon\to 0$ we
obtain the Hopf algebra \UhoE2\ given by
$$
\eqalign{
&
[J',P_x]=P_y\quad,\quad
[J',P_y]=-P_x\quad,\quad
[P_x,P_y]=\hat\Xi\quad,\quad
[\hat\Xi, {\rm all}]=0\quad;
\cr
&\Delta J'=J'\otimes 1+1\otimes J'+\rhob(\hat\Xi\otimes P_x-P_x\otimes 
\hat\Xi)\quad ,
\cr 
&\Delta P_{(x,y)}=P_{(x,y)}\otimes 1+1\otimes P_{(x,y)}\quad,\quad
\Delta\hat\Xi=\hat\Xi\otimes 1+1\otimes\hat\Xi\quad;
\cr
&S[(J',\hat\Xi,P_x,P_y)]=-(J',\hat\Xi,P_x,P_y)\quad,\quad
\epsilon[(J',\hat\Xi,P_x,P_y)]=0\quad.
\cr}
\eqn\rotexnew
$$

This algebra has a cocycle extension structure. To show this, 
we make the non-linear change
$$
J=J'+\rhob\Xi P_x\quad.\eqn\changeii
$$
This modifies only
$$
[J,P_y]=-P_x+\rhob\hat\Xi^2
\quad,\quad
\Delta J=J\otimes 1+1\otimes J+2\rhob\hat\Xi\otimes P_x\quad,\quad
S(J)=-J+2\rhob\hat\Xi P_x\quad.\eqn\rotexnewii
$$
If $\A$ is taken as 
the Hopf subalgebra generated by $\hat\Xi$ and $\H$ is the 
undeformed Euclidean Hopf algebra 
$\EUC$,
the algebra \rotexnew, \rotexnewii\ is obtained as the
right-left cocycle bicrossproduct
with $\alpha$ and $\psi$ trivial and $\beta$ and $\xi$ defined by
$$
\beta(J)=1\otimes J +2\rhob\hat\Xi\otimes P_x\quad,\quad
\xi(P_x,P_y)={1\over 2}\hat\Xi\quad,\quad
\xi(J,P_y)={\rhob\over 2}\hat\Xi^2\quad.
\eqn\newbicross
$$

\medskip
{\bf (a')}
Let us go back to the case (a) above.
If we make the redefinitions $\bar N=iN\,,\,
A=i{(X_1-iX_2)\over\sqrt 2}$, $A^{+}=i{(X_1+iX_2)\over\sqrt 2}$, the
\tUEO\ algebra in the basis \heis\ takes the form \refmark{\CGSTa,\BAM}
$$\eqalign{
&
[\bar N,A]=-A\quad,\quad
[\bar N,A^{+}]=A^{+}\quad,\quad 
[A,A^{+}]=-i{\sinh 2\rhoa\Xi\over 2\rhoa}\quad,\quad 
[\Xi,{\rm all}]=0\quad;
\cr
&
\Delta \bar N=\bar N\otimes 1+1\otimes\bar N\quad,
\quad
\Delta A=\exp(-\rhoa\Xi)\otimes A+ A\otimes \exp(\rhoa\Xi)\quad,
\cr
&
\Delta A^{+}=\exp(-\rhoa\Xi)\otimes A^{+}+ A^{+}\otimes \exp(\rhoa\Xi)
\quad,\quad
\Delta\Xi=\Xi\otimes 1+1\otimes\Xi\quad;
\cr
&
S(\bar N)=-\bar N\quad,\quad
S(A)=-A\quad,\quad S(A^{+})=-A^{+}\quad,\quad 
\epsilon(\bar N,A,A^+,\Xi)=0\quad.
\cr}
\eqn\heisa
$$

{\bf (b')}
Similarly, the redefinitions 
$\hat J=iJ\,,\,
\hat A=i{P_x-iP_y\over\sqrt 2}
\,,\,
\hat A^+=i{P_x+iP_y\over\sqrt 2}$ 
take the algebra \UhoE2\ in the basis \rotexnewii\ to the form
$$\eqalign{
&
[\hat J,\hat A]=-\hat A\quad,\quad
[\hat J,A^{+}]=\hat A^{+}\quad,\quad
[A,A^{+}]=-i\hat\Xi\quad,\quad
[\Xi,{\rm all}]=0\quad;\cr
&
\Delta \hat J=\hat J\otimes 1+1\otimes \hat J+{\sqrt 2\rhob}
\hat\Xi\otimes (\hat A +\hat A^+)\quad,
\cr
&
\Delta \hat A=1\otimes \hat A+ \hat A\otimes 1\quad,\quad
\Delta \hat A^+=1\otimes \hat A^+ + \hat A^+\otimes 1 \quad,\quad
\cr
&
S(\hat J)=-\hat J +{\sqrt 2\rhob}
\hat\Xi(\hat A+\hat A^+)\quad,\quad
S(\hat A)=-\hat A\quad,\quad 
S(\hat A^{+})=-\hat A^{+}\quad,
\cr
&
\epsilon(\hat J,\hat A,\hat A^+,\hat\Xi)=0\quad.
\cr}
\eqn\heisab$$

The algebras \heisa\ \refmark{\FIRI,\CGSTa} and \heisab\ are
a deformation of the four-generator oscillator algebra which is recovered
in the limits
$\rhoa\to 0$ [\heisa], $\rhob\to 0$ [\heisab].
Eqs. \heisa\ or \heisab\ do not, however \refmark{\CGSTa}, define the algebra
of the $q$-oscillator \refmark{\MAC,\LBIE,\AC}. The oscillator algebra 
may be obtained
by contraction using the finite-dimensional representations of $su_q(2)$
\refmark{\KUL}.
To derive it directly, without resorting to the $su_q(2)$ representations, 
consider the {\it four} generators algebra 
\su$\times$\U1\ with $[J,J_{\pm}]=\pm J_{\pm}
\;,\; [J_+,J_-]=[2J]_q\equiv\sinh 2zJ/\sinh z\;,\;[\Xi,{\rm all}]=0$. 
Now, we perform the redefinitions
\foot{
Notice that, were it not by the $q$-bracket
$[x]_q=(q^x-q^{-x})/(q-q^{-1})$, these redefinitions would be
equivalent to those in \newrescal; this exhibits once more the 
non-commutative nature of many contraction/deformation diagrams.}
$J_+=[2/\epsilon]_q^{1/2}\tilde a^+\;,\;J_-=[2/\epsilon]_q^{1/2}\tilde a\;,\;
J=N-\Xi/\epsilon\,$; this means that $\tilde a\, ,\,\tilde a^+$ and $N$ are
independent generators. Assuming $q$ real, $z>0\,,$ the contraction leads to 
$[N,\tilde a]=-\tilde a\;,\;
[N,\tilde a^{+}]=\tilde a^{+}\;,\;
[\tilde a,\tilde a^+]=q^{-2N}$; the familiar $q$-commutator relations
$[N,a]=-a\,,\,[N,a^+]=a^+\,,\,[a,a^+]_q=q^{-N}$ 
follow for $a=q^{N/2}\tilde a\;,\;a^+=\tilde a^+q^{N/2}$.
\foot{
The above oscillator algebra, where $N$ is treated as an independent generator,
has a non-trivial central element, $z=q^{-N+1}([N]_q-a^+a)$ and many
irreducible representations 
(for $0<q<1$)\refmark{\KUL} unequivalent to the Fock space ones with
vacuum state and number operator $N$, $[N]_q=a^+a,$ for which $z=0$.
}
The coproduct in \su$\times$\U1, however, does not have a limit and this 
explains why the {\it Hopf} structure for the $q$-oscillator 
(as defined by {\it these} relations) is lost 
(for recent references on this point, see \refmark{\OS,\QV}).

\chapter{The dual case: structure of the deformed Hopf group algebras}

The previous deformed algebras
may be dualized making use of the bicrossproduct
construction. 
The dual of a bicrossproduct Hopf algebra is
also a bicrossproduct Hopf algebra; thus, if $\H$ and $\A$ are Hopf algebras 
from which the bicrossproduct $\H\bic\A$ is constructed, then their duals  
$\AA$ and $\HH$ lead to the dual bicrossproduct $\AA\LR\HH$. 
This dualization
will exhibit the `group-like' (rather than `algebra-like') aspects
of the deformation. In fact, this procedure of obtaining the duals of certain 
deformed Hopf algebras is quite an efficient one, since the 
construction often embeds the non-commuting properties in some of the 
$(\alpha,\beta,\xi,\psi)$ operations, while the original algebras $\H$ and $\A$ 
are often undeformed or easy to dualize. 
We may even follow a step by step procedure.
\medskip
{\bf (1)}
The case of $\FqE2$ has been discussed in \refmark{\SH}, and will not be 
repeated here.

{\bf (2)}
Consider now \UEO\ [sec. 2(2)] which has a bicrossproduct structure according 
to sec. 3(2).
We now show that the dual algebra $\FoE2$ \refmark{\BCGST,\PM,\allowbreak
\SOB} is easily 
recovered by looking at its bicrossproduct structure.
We take for $\HH$ the dual algebra $\FoE2$ of ${\cal U}_\omega({\cal T}r(2))$ 
defined by
$$
\eqalign{
&
\Delta x=x\otimes 1+1\otimes x\quad,\quad 
\Delta y=y\otimes 1+1\otimes y\quad,\quad
[x,y]=-\omega y\quad,
\cr
&
S(x,y)=(-x,-y)\quad,\quad
\epsilon(x,y)=0\quad,\quad
x,y\in\HH\quad,
\cr}
\eqn\insertx
$$
and $\AA$ is generated by $\varphi$ with
$$
\Delta \varphi=\varphi\otimes 1+1\otimes\varphi\quad.
\eqn\insertxi
$$
The duals $\bbeta$ and $\balpha$ of $\alpha$ and $\beta$ are found to be
$$
\eqalign{
\bbeta(x)=x\otimes\cos\varphi+y\otimes\sin\varphi\quad,\quad
&
\bbeta(y)=-x\otimes\sin\varphi+y\otimes \cos\varphi\quad;
\cr
\balpha(x\otimes\varphi)\equiv x\RRR\varphi=-\omega\sin\varphi\quad,\quad
&
\balpha(y\otimes\varphi)\equiv y\RRR\varphi=\omega(1-\cos\varphi)\quad.
\cr}
\eqn\insertxii
$$
The compatibility conditions \aapix-\aapxiii\ are satisfied, and
\aapxiv, \aapxv, \aapxvi\ and \aapxvii\ 
determine the Hopf structure of $\FoE2$,
$$
\eqalign{
&
[x,y]=-\omega y\quad,\quad
[x,\varphi]=-\omega\sin\varphi\quad,\quad
[y,\varphi]=\omega(1-\cos\varphi)\quad,
\cr
&
\Delta \varphi=\varphi\otimes 1 +1\otimes\varphi\quad,\quad
\Delta x=1\otimes x+x\otimes\cos\varphi+y\otimes\sin\varphi\quad,
\cr
&
\Delta y=1\otimes y-x\otimes\sin\varphi+y\otimes\cos\varphi\quad,\quad
\epsilon(\varphi;x,y)=0\quad,
\cr
&
S(x)=-\cos\varphi x+\sin\varphi y\quad,\quad
S(y)=-\sin\varphi x -\cos\varphi y\quad.
\cr}
\eqn\insertxiii
$$

{\bf (3)}
A discussion of the Galilei case will be presented elsewhere.

{\bf (4)}
Consider now the case of the {\it deformed 
Heisenberg-Weyl `group'}\allowbreak\ $\FraHW$, 
(see \refmark{\HLR}) dual of the algebra ${\cal U}_\rhoa(HW)$ 
as given in \hei\ (\ie, \VIii\ for $\omm=2\rhoa$).
It was shown in sec. 3 (4b) that 
${\cal U}_\rho(HW)$ could be obtained as the cocycle bicrossproduct
\refmark{\MB,\MAJSOB} (Appendix B) 
${\cal U}_\rho(HW)=\H\bic\A$ of the undeformed algebras
$\H={\cal U}({\cal T}r(2))$ and $\A={\cal U}(u(1))$ by using the non-trivial 
$\beta$ and $\xi$ given \VIiv, \VIv.
Thus, the deformed Heisenberg-Weyl group algebra $\FraHW$ 
may be found as the cocycle 
bicrossproduct of $\AA={\rm Tr}(2)$ and $\HH=U(1)$ using the duals 
$\balpha:\HH\otimes\AA\to\AA$ and $\bpsi:\HH\to\AA\otimes\AA$ of $\beta$ and 
$\xi$ respectively.
Using $(y_1,y_2;\chi)$ for the parameters of ${\rm Tr}(2)$ and $U(1)$, 
$<Y_i,y_j>=\delta_{ij}\ ,\ <\Xi,\chi>=1$ the dualization of $\beta$ immediately 
leads to
$$
\balpha(\chi,y_i)\equiv\chi\RRR y_i=-2\rho y_i\quad {\rm or} \quad
[\chi,y_i]=-2\rho y_i\quad,
\eqn\inserti
$$
$i=1,2$.
Let us now dualize $\xi={1-\exp(-4\rho\Xi)\over 8\rho}$.
What was really needed in \VIv\ to compute $[Y_1,Y_2]$ was the difference 
$\xi(Y_1,Y_2)-\xi(Y_2,Y_1)\,;$
the ambiguity in $\xi(Y_1,Y_2)$ is related to the coboundary ambiguity.
A suitable election produces
$$
\bpsi(\chi)={1\over 2}(y_1\otimes y_2-y_2\otimes y_1)\quad,
\eqn\insertii
$$
from which $\Delta(\chi)$ is easily found using \aapxxxx\ since $\bbeta$ is 
trivial $(\alpha$ is trivial).
In all, $\FraHW$ is determined by
$$
\eqalign{
[y_i,y_j]=0\quad,\quad  
&
[\chi,y_i]=-2\rho y_i\quad,
\cr
\Delta y_i=y_i\otimes 1+1\otimes y_i\quad,\quad
&
\Delta\chi=\chi\otimes 1+1\otimes\chi+{1\over 2}(y_1\otimes y_2-y_2\otimes 
y_1)\quad;
\cr
S(y_i)=-y_i\quad,\quad
&
S(\chi)=-\chi\quad,\quad
\epsilon(y_i,\chi)=0\quad.
\cr}
\eqn\hwrho
$$
The coproduct mimics the familiar ${\rm HW}$ 
group law, and the non-commutativity is 
just reflected in the non-zero $[\chi,y_i]$ commutator.
\medskip
{\bf (a)}
{\it Extended Euclidean group $\FraE2$.}\quad
The dual $\FraE2$ 
of the algebra (a) given by eqs. \hei\ 
(and \heis) is generated by the elements $(y_1,y_2,\chi,\varphi)$ 
$(<Y_i,y_j>=\delta_{ij}\,,\,<N,\varphi>=1\,,\,<\Xi,\chi>=1)$
for which
$$
\eqalign{
&
[y_1,\varphi]=[y_2,\varphi]=[\chi,\varphi]=[y_1,y_2]=0\quad,
\cr
&
[\chi,y_1]=-2\rhoa y_1\quad,\quad [\chi,y_2]=-2\rhoa y_2\quad;
\cr
&
\Delta \varphi=\varphi\otimes 1+1\otimes\varphi\quad,\quad
\Delta y_1=1\otimes y_1+y_1\otimes \cos\varphi+y_2\otimes\sin\varphi\quad,
\cr
&
\Delta y_2=1\otimes y_2+y_2\otimes \cos\varphi-y_1\otimes\sin\varphi\quad,
\cr
&
\Delta\chi=1\otimes\chi+\chi\otimes 1
\cr
&
{\mskip 40mu}
+{1\over 2} [y_1\otimes \cos\varphi y_2+ y_2\otimes \sin\varphi y_2 
- y_2\otimes \cos\varphi y_1 + y_1\otimes \sin\varphi y_1]
\,;
\cr
&
S(\varphi)=-\varphi\quad,\quad 
S(y_1)=-\cos\varphi y_1+\sin\varphi y_2\;,\;
S(y_2)=
-\sin\varphi y_1-\cos\varphi y_2
\;,\cr
&
S(\chi)=-\chi\quad;\quad
\epsilon(\varphi,y_1,y_2,\chi)=0\quad.
\cr}
\eqn\duali
$$
It is not difficult to check directly that $\FraE2$ [\duali] 
is a Hopf algebra;
we shall now obtain \duali\ by dualization in two different ways.
For the dual in the basis \heisa, see \refmark{\BAM}.
\medskip
\noindent
a1) $\FraE2$ 
is the bicrossproduct ${\rm Fun}U(1)\LR\FraHW$, where $U(1)$ is 
generated by $\varphi$ and $\FraHW$ is given in \hwrho.
To see this, it is sufficient to dualize the right action $\LLL$ (eq. \acc),
$Y_1 \LLL N=-Y_2\,,\, Y_2\LLL N=Y_1\,,\, \Xi\LLL N=0$ to obtain
$$
\eqalign{
\bbeta(y_1)=y_1\otimes \cos\varphi+y_2\otimes\sin\varphi\quad,
&
\quad
\bbeta(y_2)=-y_1\otimes\sin\varphi+y_2\otimes \cos\varphi\quad,
\cr
\bbeta(\chi)=
&
1\otimes\chi\quad,
\cr}
\eqn\insertiii
$$
for which the coproducts and antipodes in \duali\ 
are obtained from \aapxv\ and \aapxvii. 
Clearly $[y_i,\varphi]=0=[\chi,\varphi]$ since $\balpha\ (\RRR)$ 
is dual to $\beta$,
which is trivial.
\medskip
\noindent
a2)
$\FraE2$ has also a cocycle bicrossproduct structure.
To see this, 
we take $\HH$ as the Hopf algebra 
generated by $\chi$ with primitive coproduct 
and $\AA$ as the (undeformed) Euclidean group Hopf algebra 
of generators $(y_1,y_2,\varphi)$ with $\Delta y_i\,,\Delta\varphi\ ,
\allowbreak\, 
S(y_i,\varphi)$ and 
$\epsilon(y_i,\varphi)$ as in \duali.
Then, since $\alpha$ and $\psi$ were trivial in sec. 4 a2), 
$\bbeta$ and $\bxi$ are trivial ($\bbeta(a)=a\otimes 
1_\AA\,,\,\bxi(a,b)=\epsilon(a)\epsilon(b)1_\AA$) and 
$\bpsi:\HH\to\AA\otimes\AA\,,\,\balpha:\HH\otimes\AA\to\AA$ may be found from 
\defcoci\ to be
$$
\eqalign{
&\bpsi(\chi)= {1\over 2}[
y_1\otimes \cos\varphi y_2 + y_2\otimes \sin\varphi y_2 
-y_2\otimes\cos\varphi y_1 +y_1\otimes\sin\varphi y_1] \quad,
\cr
&\chi\RRR y_1=-2\rhoa y_1\quad,\quad
\chi\RRR y_2=-2\rhoa y_2\quad,\quad
\chi\RRR \varphi= 0\quad.
\cr
}
\eqn\dualii
$$
The relations \aapxx-\aapxxviiia\ are fulfilled
\foot{
The only non-trivial properties are \aapxxiv\ and \aapxxvi. The first one is 
the dual cocycle condition, verified because the dual cocycle $\bpsi$ is 
the undeformed one, and the second one is due to the compatibility between
the coproduct and the commutators.
}
and the cocycle 
bicrossproduct structure of $\FraE2$ follows from \aapxxix\ (which for 
$\bxi$ trivial and $\hh$ with primitive coproduct leads to 
$[\hh,\aa]=\hh\RRR\aa$) and \aapxxxx.

\medskip
{\bf (b)}
{\it Extended Euclidean group $\FrbE2$}.\quad
This is the dual $\FrbE2$ of
the Hopf algebra \UhoE2\ (see eqs. \rotexnew\ and 
\rotexnewii).
It is generated by the elements $(x,y,\varphi,\hat\chi)$ 
$(<P_x,x>=<P_y,y>=<J,\varphi>=<\hat\Xi,\hat\chi>=1)$
with relations
$$
\eqalign{
&
[x,\varphi]=[y,\varphi]=[\hat\chi,\varphi]=[x,y]=0\quad,
\cr
&
[\hat\chi,x]=-2\rhob \sin\varphi
\quad,\quad 
[\hat\chi,y]=2\rhob (1-\cos\varphi)\quad;
\cr
&
\Delta \varphi=\varphi\otimes 1+1\otimes\varphi\quad,\quad
\Delta x=1\otimes x+x\otimes \cos\varphi+y\otimes\sin\varphi\quad,
\cr
&
\Delta y=1\otimes y+y\otimes \cos\varphi-x\otimes\sin\varphi\quad,
\cr
&
\Delta\hat\chi=1\otimes\hat\chi+\hat\chi\otimes 1 + {1\over 2}[
x\otimes \cos\varphi y +y\otimes \sin\varphi y 
-y\otimes \cos\varphi x+x\otimes \sin\varphi x] 
\;,
\cr
&
S(\varphi)=-\varphi\quad,\quad
S(x)=-\cos\varphi x+\sin\varphi y\quad,\quad
S(y)=-\cos\varphi y-\sin\varphi x\quad,
\cr
&
S(\chi)=-\chi\quad;\quad
\epsilon(\varphi,x,y,\chi)=0\quad,
\cr}
\eqn\dualx
$$
which define a Hopf algebra as it may be checked.
Now, we take $\HH$ as the Hopf group algebra generated by $\hat\chi$ and
$\AA$ as 
the dual undeformed Euclidean group Hopf algebra ${\rm Fun}(E(2))$
(eqs. \insertxiii\ for $\omega=0$) 
of generators $(x,y,\varphi)$. 
If we now define $\bbeta\,,\,\bxi$ to be trivial plus
$$
\eqalign{
&\bpsi(\hat\chi)=
{1\over 2}[x\otimes\cos\varphi y+y\otimes\sin\varphi y-
y\otimes\cos\varphi x +x\otimes\sin\varphi x] 
\quad,
\cr
&\hat\chi\RRR x=-2\rhob\sin\varphi \quad,\quad
\hat\chi\RRR y=2\rhob(1-\cos\varphi)\quad,\quad
\hat\chi\RRR \varphi= 0\quad,
\cr}
\eqn\dualxi
$$
from $\xi$ and $\beta$ in eq. \newbicross, 
the Hopf algebra \dualx\ is recovered using \aapxxix\ and \aapxxxx,
which exhibits the cocycle 
bicrossproduct structure of \dualx.
Due to the commutators $[\chi,x]\,,\,[\chi,y]$, there is no Hopf 
$\FrbHW$ subalgebra here and no bicrossproduct structure in contrast with 
the previous a1) case.
\chapter{Differential calculus on the Euclidean and Galilean planes}

We shall now introduce a covariant differential calculus \refmark{\WORO} (see 
Appendix A) on the different homogeneous spaces which can be constructed.
Clearly, to have a proper action on the `homogeneous' part, a bicrossproduct 
structure is needed. Let us consider now a few different cases.

{\bf (1)}
Due to the lack of a bicrossproduct structure, the inhomogeneous part of the
\UE2\ algebra $(P_1,P_2)$ does not constitute a Hopf 
subalgebra, and the 
construction of the space algebra as the dual of $(P_1,P_2)$
cannot be performed.

{\bf(2)}
The Euclidean plane $E^2_\omega$
is introduced as the dual ($<P_i,x_j>=\delta_{ij}$)
of the translation Hopf subalgebra 
${\cal U}_\omega({\cal T}r(2))$ 
of \UEO\  
generated by $P_i$ (eq. \euclid).
Since ${\cal U}_\omega({\cal T}r(2))$ 
is commutative but not cocommutative, we obtain (eqs. \insertx)
$$
\Delta x=x\otimes 1+1\otimes x\quad,\quad
\Delta y=y\otimes 1+1\otimes y\quad;\quad
[x,y]=-\omega y\quad;
\eqn\plano
$$
for the $E^2_{\omega}$-plane algebra 
associated with \UEO.
Let us construct a bicovariant differential calculus on
$E^2_{\omega}$ which is consistent (\ie\ covariant) under the action of $J$.
The (left) action of $J$ on $E^2_{\omega}$ is defined by duality,
$<P_x\acti J,x>=<P_x,J\act x>$ etc., 
from which follows that
$$
J\act x=y\quad,\quad J\act y=-x\quad.
\eqn\covariant
$$
To define a first order ($J$-)covariant differential calculus 
we have to determine all commutators $[x_i,dx_j]$ in a way which is 
consistent with the action \covariant\ (which for instance, implies 
$J\act xdy=(J_{(1)}\act x) d(J_{(2)}\act y)$) and with the Jacobi identity.
Although it is not difficult to check that the set of covariance equations 
(like
$J\act(x_idx_j)-J\act(dx_jx_i)=J\act[x_i,dx_j])$ has a {\it unique} solution
given by
$$[x,y]=-\omega y\quad,\quad[x,dx]=0=[x,dy]\quad,\quad[y,dx]=\omega dy
\quad,\quad[y,dy]=-\omega dx\quad,\eqn\solut$$ 
the above commutators do not satisfy the Jacobi identity and thus fail to
provide a consistent differential calculus.
This situation is not new, and has already appeared for the differential 
calculus on other spacetime algebras 
\refmark{\SIT,\APB}.
We now show that the solution proposed there, and which involves an 
enlargement of the algebra which has been found to be associated with a Hopf
algebra cocycle extension \refmark{\APB}, also applies here.
We stress that this problem is associated to the deformed character of 
\euclid\ as expressed by $\omega$, being of course 
absent for the undeformed Euclidean Hopf algebra ${\cal U}({\cal E}_2)$.

Consider the trivial extension \UEO$\times$\U1\ mentioned in sec. 4, obtained 
by adding the primitive Hopf algebra generated by $\Xi$ 
to \UEO.
The previous procedure applied to $(P_x,P_y,\Xi)$ leads now 
to an enlarged Euclidean algebra ${\tilde E}_{\omega}$
generated by $(x,y,\chi)$ $(<\Xi,\chi>=1)$ and to the additional
relations
$$[\chi,x]=0=[\chi,y]\quad,\quad\Delta\chi=\chi\otimes 1+1\otimes\chi\quad,
\quad J\act\chi=0\quad.\eqn\defchi$$
Proceeding as before, we find that there is a unique solution 
for the rotation covariant differential calculus on the above enlarged 
Euclidean `space' specified by (cf. \solut)
$$\eqalign{&[x,y]=-\omega y\ ,\ [x,\chi]=0\ ,\ 
[x,dx]=\omega d\chi\ ,\ [x,dy]=0\ ,\  [x,d\chi]=\omega dx\ ,\cr
&[y,\chi]=0\quad,\quad[y,dx]=\omega dy\quad,\quad [y,dy]=-\omega(dx-d\chi)
\quad,\quad [y,d\chi]=\omega dy\quad,\cr
&[\chi,dx]=\omega dx\quad,\quad[\chi,dy]=\omega dy\quad,\quad 
[\chi,d\chi]=\omega d\chi\quad,\cr}\eqn\solution$$
and satisfying Jacobi identity.

{\bf (3)}
We define the two-dimensional Galilean plane $\galp$
as the dual ($<X,x>=1=<X_t,t>\,,\,<X,t>=0=<X_t,x>$)
of $\trrr$.
The commutativity (non-cocommutativity) of $\trrr$ implies the relations
$$
[x,t]=2\tilde\omega x\quad;\quad
\Delta x=x\otimes 1+1\otimes x\quad,\quad
\Delta t=t\otimes 1+1\otimes t\quad,
\eqn\galplane
$$
for the $\galp$ algebra. 
Following the same pattern of case (2)
we construct a bicovariant differential calculus 
(covariant under the action of the `boost' $V$) that satisfies Leibniz's 
rule and Jacobi identity.
The (left) action of $V$ on $\galp$ is given by
$$
V\act x=-t\quad,\quad V\act t=0\quad.
\eqn\covargal
$$
Using \covargal, 
we find that the covariance requirement implies the system 
of equations
$$
\eqalign{
V\act [x,dx]=-[t,dx]-[x,dt]\quad ,&\quad 
V\act[t,dx]=-[t,dt]+2\tilde\omega dt
\cr
V\act[x,dt]=-[t,dt]\quad,&\quad
V\act[t,dt]=0\quad.\cr}
\eqn\syst
$$
The unique solution linear in $dx, dt$ that satisfies \syst, Leibniz's rule 
and Jacobi identity is
\foot
{Even if 
there is no deformation ($\tilde\omega=0$) there exists a non-trivial 
solution (see \refmark{\APB}) given by 
$[x,dx]=\mu dt$
and all other commutators equal to zero.
}
$$
[x,dx]=0\quad,\quad 
[x,dt]=\tilde\omega dx\quad,\quad 
[t,dx]=-\tilde\omega dx\quad,\quad
[t,dt]=\tilde\omega dt\quad.
\eqn\galsolution
$$
Thus, this case
is different from the Euclidean case $E^2_\omega$.
On $\galp$ there is a covariant differential calculus without any 
additional one-form.
\foot{
For the differential calculus on the deformed Newtonian spacetime associated 
with the $(1+3)$ version of the deformed Galilei algebra 
${\cal G}_{\tilde\kappa}$ see \refmark{\APB}.
}

{\bf (a)}
It was seen (eq. \solution) that to define a 
$J$-covariant differential calculus on $E^2_\omega$ 
it was necessary to enlarge it to ${\tilde E}_\omega$. 
Let us now show that two $N$-covariant calculi may be similarly constructed 
on $\FraHW$ (eqs. \hwrho) as the dual of the
${\cal U}_\rhoa(HW)$ subalgebra of \tUEO\ (sec. 4(a)).
The left action $\act$ of $N$ on $(y_1,y_2,\chi)$ is
obtained from \acc\ and given by
$$
N\act y_1=y_2\quad,\quad N\act y_2=-y_1\quad,\quad N\act\chi=0\quad.
\eqn\heiac
$$
Proceeding as before, we find the commutators 
$$
\eqalign{
&
[y_i,y_j]=0\quad,\quad 
[y_i,\chi]=2\rhoa y_i\quad,\quad
[y_i,dy_j]=0\quad,\cr
&[\chi, dy_i]=\lambda dy_i\quad,\quad 
[y_i,d\chi]=(\lambda+2\rhoa)dy_i
\quad,\quad 
[\chi,d\chi]=\mu d\chi\quad.\cr}\eqn\ecui$$
The Jacobi identity requires $\lambda=-2\rhoa\ {or}\ \lambda-\mu=-2\rhoa$.
The bicovariance requirement now determines two bicovariant 
differential calculi over $\FraHW$ (on the  
$E^2_\omega$ plane the coproduct of the generators was primitive, hence
the differentials are bi-invariant by \difoo\ and the bicovariance is trivial).
We first find, using \hwrho\ and \difoo, 
$$\eqalign{\Delta_L dy_i=1\otimes dy_i\quad,
&\quad\Delta_R dy_i=dy_i\otimes 1
\quad,\cr
\Delta_L d\chi=1\otimes d\chi+{1\over 2}
&(y_1\otimes dy_2-y_2\otimes dy_1)\quad,\cr 
\Delta_R d\chi=d\chi\otimes 1+{1\over 2}
&(dy_1\otimes y_2-dy_2\otimes y_1)\quad;\cr}
\eqn\ecudos$$
it is easy to show that the coactions \ecudos\ satisfy \difi. 
If we use now \difo\ to calculate $\Delta_L[\chi,d\chi]=\mu\Delta_L d\chi$ 
we find $\mu=2\lambda$; the same condition is obtained 
using $\Delta_R$. Then, \ecui\ leads to $(\lambda=-2\rhoa,\mu=-4\rhoa)$
$$
\eqalign{
&
[y_i,y_j]=0\quad,\quad 
[y_i,\chi]=2\rhoa y_i\quad,\quad
[y_i,dy_j]=0\quad,
\cr
&[\chi, dy_i]=-2\rhoa dy_i\quad,\quad [y_i,d\chi]=0
\quad,\quad 
[\chi,d\chi]=-4\rhoa d\chi
\cr}
\eqn\ecutres
$$
and $(\lambda=2\rhoa,\mu=4\rhoa)$
$$
\eqalign{
&
[y_i,y_j]=0\quad,\quad 
[y_i,\chi]=2\rhoa y_i\quad,\quad
[y_i,dy_j]=0\quad,
\cr
&[\chi, dy_i]=2\rhoa dy_i\quad,\quad 
[y_i,d\chi]=4\rhoa dy_i\quad,\quad 
[\chi,d\chi]=4\rhoa d\chi\quad.
\cr}
\eqn\ecucuatro
$$
Since \difii\ is satisfied, eqs. \ecutres, \ecucuatro\ 
determine two first order $N$-covariant differential calculi over \tEtO.

\ack
This paper has been partially supported by the CICYT grant AEN93-187.
One of us (JCPB) wishes to acknowledge a FPI grant from the Spanish Ministry 
of Education and Science and the CSIC.
Both authors wish to thank M. del Olmo for very helpful discussions.
\vskip 24 pt                                                             
\noindent
{\twelvebf Appendix A: Bicovariant differential calculus}
\ABC{A}

Let $\HH$ be a Hopf algebra and let $\Delta$ and $\epsilon$ be its coproduct 
and counit. A first order bicovariant differential
calculus over $\HH$ is defined \refmark{\WORO} 
by a pair $(\Gamma,d)$ where $d:\HH\to\Gamma$ is a linear mapping satisfying 
Leibniz's rule and $\Gamma$ is a bicovariant $\HH$-bimodule
$(\Gamma,\Delta_L,\Delta_R)$ \ie, the linear mappings $\Delta_L:
\Gamma\rightarrow \HH\otimes\Gamma\;,\; 
\Delta_R:\Gamma\rightarrow\Gamma\otimes \HH$ and the exterior 
derivative $d$ satisfy 
$$\eqalign{\Delta_L(a\omega)=\Delta(a)\Delta_L(\omega)\quad, &\quad
\Delta_R(a\omega)=\Delta(a)\Delta_R(\omega)\quad,\cr
\Delta_L(\omega a)=\Delta_L(\omega)\Delta(a)\quad, &\quad 
\Delta_R(\omega a)=\Delta_R(\omega)\Delta(a)\quad,\cr}\eqn\difo$$
$$\eqalign{(\Delta\otimes id)\Delta_L=(id\otimes\Delta_L)\Delta_L\quad,&\quad
(id\otimes\Delta)\Delta_R=(\Delta_R\otimes id)\Delta_R\quad,\cr
(\epsilon\otimes id)\Delta_L=id\quad &,\quad(id\otimes\epsilon)\Delta_R=id
\quad,\cr}\eqn\difi$$
$$\eqalign{(id\otimes\Delta_R)\Delta_L=&
(\Delta_L\otimes id)\Delta_R\quad;\cr}\eqn\difii$$
$$\Delta_L d=(id\otimes d)\Delta\quad,\quad\Delta_R d=(d\otimes id)
\Delta\quad,\eqn\difoo$$
where the left (right) equations in \difi\ express that $\Gamma$ is a left
(right) $\HH$-comodule, \difii\ is the result of bicovariance (commutation of
the left and right coactions), and \difoo\ expresses the compatibility of the
exterior derivative $d$ with $\Delta$ and $\Delta_{L,R}$.
Eqs. \difo, \difi\ and \difii\ characterize $(\Gamma,\Delta_L,\Delta_R)$ as a 
bicovariant bimodule over $\HH$; the addition of \difoo\ determines a first 
order bicovariant differential calculus $(\Gamma, d)$. 
An element $\omega\in\Gamma$ is called left (right) invariant if 
$\Delta_L(\omega)=1\otimes\omega\ (\Delta_R(\omega)=\omega\otimes 1)$.
As in the undeformed (Lie) case, the basis elements of the vector space 
$\Gamma_{\rm inv}\subset\Gamma$ of the left-invariant elements generate $\Gamma$
as a left free module. 
\vskip 24pt
\noindent
{\twelvebf Appendix B: Bicrossproduct of Hopf algebras and cocycles}
\ABC{B'}

We list here for convenience the basic formulae of 
Majid's bicrossproduct and cocycle 
bicrossproduct constructions and refer to 
\refmark{\MB,\MAJSOB} (see also \refmark{\BCM}) for details.
The expressions which characterize $\H\RL\A$ (used in secs. 3,4) involve the 
mappings $\alpha:\A\otimes\H\to\A$ (right $\H$-module action), 
$\beta:\H\to\A\otimes\H$ (left $\A$-comodule coaction),
$\xi:\H\otimes\H\to\A$ (two-cocycle) and $\psi:\H\to\A\otimes\A$ 
(hence the more
detailed notation $\H^{\beta,\psi}\RL_{\alpha,\xi}\A$, see \refmark{\MB}).
Those of the dual case ( $\AA_{\balpha,\bxi}\LR^{\bbeta,\bpsi}\HH$ when all
ingredients are indicated) involve the respective dual operations; 
they were used in sec. 5.
We may think of $\H\bic\A$ as emphasizing the
`algebra-like' aspects and of $\AA\LR\HH$ 
as giving the `group-like' ones
\foot{
We use the bicrossproduct notation $\RL$ or $\RL$ rather than the (right, 
left) crossproduct ($\RIMO,\LEMO$) or the (left, right) cross coproduct 
($\LECO,\RICO$) even if the coactions $\beta,\bbeta$ or the actions 
$\alpha,\balpha$ are trivial, and omit explicit reference to them (or to
$\xi,\psi$ etc.)} 
(in the undeformed case they correspond, respectively, 
to the cocommutative Hopf algebra 
constructed on the enveloping algebra ${\cal U}({\cal G})$ of ${\cal G}$,
and to the Abelian Hopf algebra of 
functions Fun$(G)$ 
over a Lie group $G$ with coproduct given by the group law).
Both sets of formulae are in correspondence once $\H$, $\A$ 
$(\alpha,\beta,\xi,\psi)$ are replaced by their respective duals $\AA$, 
$\HH$, $(\bbeta,\balpha,\bpsi,\bxi)$;
thus we shall only reproduce here those for the second case.
Those useful for $\H\RL\A$ may be found in the original papers
\refmark{\MB,\MAJSOB}
(or in the Appendix of \refmark{\APB} with the same numbers 
they are referred to in the main text, 
also corresponding to the dual formulae for $\AA\LR\HH$ below).

Let $\AA$ and $\HH$ be Hopf algebras and let
\mypointbegin
$\AA$ be a left $\HH$-module algebra ($\AA\LEMO \HH$)
\mypoint
$\HH$ be a right $\AA$-comodule coalgebra ($\AA\RICO \HH$)
\ie, there exist linear mappings
$$\balpha:\HH\otimes \AA\rightarrow \AA\quad,\quad\balpha(\hh\otimes \aa)
\equiv \hh\RRR \aa
\quad,\quad \hh\in \HH,\ \aa\in \AA\quad;\eqn\aapi$$
$$\bbeta:\HH\rightarrow \HH\otimes \AA\quad,\quad \bbeta(\hh)=\hh^{(1)}\otimes 
\hh^{(2)}\quad,\quad \hh^{(1)}\in \HH,\ \hh^{(2)}\in \AA\eqn\aapii$$
such that the properties of

\noindent
a1) $\balpha$ being a left $\HH$-module action $\RRR$:
$$1_{\HH}\RRR \aa=\aa\quad,\eqn\aapiiia$$
$$\hh'\RRR (\hh\RRR \aa)=\hh'\hh\RRR \aa\quad;
\eqn\aapiii$$
a2) $\AA$ being a left $\HH$-module algebra:
$$\hh\RRR 1_{\AA}=1_{\AA}\epsilon(\hh)\quad,\quad 
\hh\RRR(\aa\bb)=(\hh_{(1)}\RRR \aa)(\hh_{(2)}\RRR \bb)\quad;
\eqn\aapiv$$
b1) $\bbeta$ being a right $\AA$-comodule coaction:
$$\hh^{(1)}\otimes \epsilon_{\AA}(\hh^{(2)})=\hh\otimes 1_{\AA}\equiv \hh\quad
[(id\otimes\epsilon)\circ\beta=id]\quad,\eqn\aapv$$
$$\hh^{(1)(1)}\otimes \hh^{(1)(2)}\otimes \hh^{(2)}=\hh^{(1)}
\otimes \hh^{(2)}_{\ (1)}\otimes \hh^{(2)}_{\ (2)}
\quad [(\bbeta\otimes id)\circ\bbeta=
(id\otimes\Delta)\circ\bbeta]\quad;\eqn\aapvi$$
b2) ${\HH}$ being a right ${\AA}$-comodule coalgebra:
$$\epsilon_{\HH}(\hh^{(1)})\hh^{(2)}=1_{\AA}\epsilon_{\HH}(\hh)\quad
[(\epsilon\otimes id)\circ \bbeta=\epsilon]\quad,\eqn\aapvii$$
$$\eqalign{& 
\hh^{(1)}_{\ (1)}\otimes \hh^{(1)}_{\ (2)}\otimes \hh^{(2)}=
\hh_{(1)}^{\ (1)}\otimes \hh_{(2)}^{\ (1)}\otimes \hh_{(1)}^{\ (2)}
\hh_{(2)}^{\ (2)}
\quad\cr 
[(\Delta\otimes id)\circ\bbeta&=
(id\otimes id\otimes m_{\AA})\circ(id\otimes\tau\otimes id)\circ
(\bbeta\otimes\bbeta)\circ\Delta
\equiv(\bbeta\hat\otimes\bbeta)\circ\Delta]\;,\cr}\eqn\aapviii$$
where $m_{\AA}$ is the multiplication in ${\AA}$ and $\tau$ is 
the twist mapping, are fulfilled.

Then, if the compatibility conditions
$$\epsilon_{\AA}(\hh\RRR \aa)=\epsilon_{\HH}(\hh)\epsilon_{\AA}(\aa)\quad,
\eqn\aapix$$
$$\Delta(\hh\RRR \aa)\equiv (\hh\RRR \aa)_{(1)}\otimes (\hh\RRR \aa)_{(2)}
=(\hh_{(1)}^{\ (1)}\RRR \aa_{(1)})\otimes \hh_{(1)}^{\ (2)}(\hh_{(2)}
\RRR \aa_{(2)})
\quad,\eqn\aapx$$
$$\bbeta(1_{\HH})\equiv 1_{\HH}^{(1)}\otimes 1_{\HH}^{(2)}=
1_{\HH}\otimes 1_{\AA}\quad,\eqn\aapxi$$
$$\bbeta(\hh\gg)\equiv (\hh\gg)^{(1)}\otimes (\hh\gg)^{(2)}
=\hh_{(1)}^{\ (1)}\gg^{(1)}\otimes \hh_{(1)}^{\ (2)}(\hh_{(2)}\RRR \gg^{(2)}) 
\quad,\eqn\aapxii$$
$$\hh_{(2)}^{\ (1)}\otimes(\hh_{(1)}\RRR \aa)\hh_{(2)}^{\ (2)}=
\hh_{(1)}^{\ (1)}\otimes \hh_{(1)}^{\ (2)}(\hh_{(2)}\RRR \aa)\quad,
\eqn\aapxiii$$
are satisfied
\foot{
If $\HH$ is cocommutative and $\AA$ commutative, condition \aapxiii\ is 
automatically satisfied.}
, 
there is a Hopf algebra structure on \refmark{\MB}
$\KK\equiv{\AA}\otimes {\HH}$ called the (left-right) bicrossproduct 
${\AA}_{\balpha}\LR^{\bbeta}{\HH}$ (${\AA}\LR {\HH}$ for short)  
defined by
$$(\aa\otimes \hh)(\bb\otimes \gg)=\aa(\hh_{(1)}\RRR \bb)\otimes 
\hh_{(2)}\gg\quad,
\quad\aa,\bb\in {\AA}\;;\;\hh,\gg\in {\HH}\quad,
\eqn\aapxiv$$
$$\Delta_{\KK}(\aa\otimes \hh)=\aa_{(1)}\otimes \hh_{(1)}^{\ (1)}
\otimes \aa_{(2)}
\hh_{(1)}^{\ (2)}\otimes \hh_{(2)}\quad,\eqn\aapxv$$
$$\epsilon_{\KK}=\epsilon_{\AA}\otimes\epsilon_{\HH}\quad,\quad
 1_{\KK}=1_{\AA}\otimes 1_{\HH}\quad,\eqn\aapxvi$$
$$S(\aa\otimes \hh)=(1_{\AA}\otimes S_{\HH}(\hh^{(1)}))(S_{\AA}(\aa\hh^{(2)})
\otimes 1_{\HH})
\quad.\eqn\aapxvii$$
In $\KK={\AA}\otimes {\HH},\;
\aa\equiv \aa\otimes 1_{\HH}$ and $\hh\equiv 1_{\AA}\otimes \hh$; thus, 
$\hh\aa=\hh_{(1)}\RRR \aa\otimes \hh_{(2)}$.
There are two cases of special interest \refmark{\MOL} (see also 
\refmark{\MB}).
When $\bbeta=I\otimes 1_\AA$ \ie\ 
$\bbeta(\hh)=\hh\otimes 1_\AA$ (trivial coaction) and $\HH$ 
is cocommutative, $\KK$ 
is the semidirect {\it product}
of Hopf algebras since then
$\Delta_\KK(\aa\otimes\hh)=
(\hh_{(1)}\otimes\aa_{(1)})\otimes(\gg_{(2)}\otimes\bb_{(2)})$.
When $\balpha$ is trivial, $\balpha=
\epsilon_\HH\otimes 1_\AA\ (\hh\RRR \aa =\aa\epsilon_\HH(\hh))$ and 
$\AA$ is commutative,
$\KK$ is the semidirect {\it coproduct} of Hopf algebras since 
$(\aa\otimes\hh)(\bb\otimes\gg)=\aa\bb\otimes\hh\gg$. 
When $\balpha$ is trivial, $\bbeta(\hh\gg)=\bbeta(\hh)\bbeta(\gg)$ 
(algebra homomorphism).

As for $\H\RL\A$, the above construction may be extended to accommodate 
cocycles \refmark{\MB, \MAJSOB}. 
Let ${\HH}$ and ${\AA}$ two Hopf algebras and $\balpha$ and $\bbeta$
as in \aapi, \aapii. Then ${\AA}$ 
is a left ${\HH}$-module cocycle algebra if \aapiiia, \aapiv\ 
are fulfilled and 
there is a linear (two-cocycle) map 
$\bxi:{\HH}\otimes {\HH}\rightarrow {\AA}$ such that
$$\bxi(\hh\otimes 1_{\HH})=1_{\AA}\epsilon(\hh)=\bxi(1_{\HH}\otimes \hh)
\quad [\bxi(1_{\HH}\otimes 1_{\HH})=1_{\AA}]\quad,\eqn\aapxx$$
$$\hh_{(1)}\RRR\bxi(\gg_{(1)}\otimes \ff_{(1)})
\bxi(\hh_{(2)}\otimes \gg_{(2)}\ff_{(2)})
=\bxi(\hh_{(1)}\otimes \gg_{(1)})\bxi(\hh_{(2)}\gg_{(2)}\otimes \ff)\;,\; 
\forall
\hh,\gg,\ff\in {\HH}\;,\eqn\aapxxi$$
(cocycle condition)
and \aapiii\ is replaced by
$$\hh_{(1)}\RRR(\gg_{(1)}\RRR \aa)\bxi(\hh_{(2)}\otimes \gg_{(2)})=
\bxi(\hh_{(1)}\otimes \gg_{(1)})((\hh_{(2)}\gg_{(2)})\RRR \aa)\quad, 
\forall \aa\in 
{\AA},\forall \hh,\gg\in {\HH}\quad,\eqn \aapxxii$$
which for $\bxi$ trivial reproduces \aapiii. Similarly, ${\HH}$ is a right 
${\AA}$-comodule coalgebra cocycle if \aapv,\ \aapvii,\ \aapviii\ 
are fulfilled, and 
there is a linear map $\bpsi:{\HH}\rightarrow {\AA}\otimes {\AA}$, 
$\bpsi(\hh)=\bpsi(\hh)^{(1)}\otimes\bpsi(\hh)^{(2)}$, such that
$$\epsilon(\bpsi(\hh)^{(1)})\bpsi(\hh)^{(2)}=1\epsilon(\hh)=
\bpsi(\hh)^{(1)}\epsilon(\bpsi(\hh)^{(2)})\ ,\ \left[(\epsilon\otimes 
id)\circ\bpsi=(id\otimes\epsilon)\circ\bpsi\right]\ ,\eqn\aapxxiii$$
$$\Delta\bpsi(\hh_{(1)})^{(1)} \bpsi(\hh_{(2)}^{\ (1)})\otimes 
\bpsi(\hh_{(1)})^{(2)}\hh_{(2)}^{\ (2)}=
\bpsi(\hh_{(1)})^{(1)}\otimes\Delta\bpsi(\hh_{(1)})^{(2)}
\bpsi(\hh_{(2)}),\ \forall \hh\in {\HH}\quad,\eqn\aapxxiv$$
(dual cocycle condition) and \aapvi\ is replaced by
$$\eqalign{(1\otimes\bpsi(\hh_{(1)}))((\bbeta\otimes id)\circ\bbeta(\hh_{(2)}))
&=\hh_{(1)}^{\ (1)}\otimes \Delta \hh_{(1)}^{\ (2)}\bpsi(\hh_{(2)})\cr
&=((id \otimes \Delta)\bbeta(\hh_{(1)}))(1\otimes\bpsi(\hh_{(2)}))
\quad.\cr}\eqn\aapxxv$$

Then, if the compatibility conditions \aapix,\ \aapxi,\ \aapxiii\ and
$$\Delta(\hh_{(1)}\RRR \aa)\bpsi(\hh_{(2)})=\bpsi(\hh_{(1)})
[\hh_{(2)}^{\ (1)}\RRR 
\aa_{(1)}\otimes \hh_{(2)}^{\ (2)}
(\hh_{(3)}\RRR \aa_{(2)})]\quad, \eqn\aapxxvi$$
$$(1\otimes\bxi(\hh_{(1)}\otimes \gg_{(1)}))\bbeta(\hh_{(2)}\gg_{(2)})=
\hh_{(1)}^{\ (1)}\gg_{(1)}^{\ (1)}\otimes \hh_{(1)}^{\ (2)}(\hh_{(2)}
\RRR \gg_{(1)}^{\ (2)})\bxi(\hh_{(3)}\otimes \gg_{(2)})
\quad,\eqn\aapxxvii$$
(which replace \aapx
\foot{
With $\bxi(\hh_{(1)}\otimes \gg_{(1)})\bxi^{-1}(\hh_{(2)}\otimes 
\gg_{(2)})=\epsilon(\hh)\epsilon(\gg)$ (convolution invertible \refmark{\MB}), eq. 
\apxxvii\ gives $\beta(\hh\gg)=
\hh_{(2)}^{\ (1)}\gg_{(2)}^{\ (1)}\otimes 
\bxi^{-1}(\hh_{(1)}\otimes \gg_{(1)})
\hh_{(2)}^{\ (2)}(\hh_{(3)}
\RRR \gg_{(2)}^{\ (2)})\bxi(\hh_{(4)}\otimes \gg_{(3)})
$.
If $\HH$ is Abelian, as is always the case in the cocycle bicrossproduct  
structures in the main text, this formula reduces to \aapxii.}
, \aapxii), together with
$$\eqalign{\Delta\bxi(\hh_{(1)}\otimes \gg_{(1)})\bpsi(\hh_{(2)}\gg_{(2)})
&=\bpsi(\hh_{(1)})\left[(\hh_{(2)}^{\ (1)}\RRR\bpsi(\gg_{(1)})^{(1)})
\bxi(\hh_{(4)}^{\ (1)}\otimes \gg_{(2)}^{\ (1)})\otimes\right.\cr
&\left. \hh_{(2)}^{\ (2)}(\hh_{(3)}\RRR\bpsi(\gg_{(1)})^{(2)})\hh_{(4)}^{\ (2)}
(\hh_{(5)}\RRR \gg_{(2)}^{\ (2)})\bxi(\hh_{(6)}\otimes \gg_{(3)})\right]
\;,\cr}\eqn\aapxxviii$$
$$\eqalign{\epsilon(\bxi(\hh\otimes \gg))&=\epsilon(\hh)\epsilon(\gg)
\quad,\quad \bpsi(1_{\HH})=1_{\AA}\otimes 1_{\AA}\cr}\eqn\aapxxviiia$$
hold, $({\HH},{\AA},\balpha,\bbeta,\bxi,\bpsi)$ 
determine a cocycle left-right bicrossproduct bialgebra \break 
${\AA} _{\bxi}\LR^{\bpsi} {\HH}$. In it, the counit and 
unit are defined by \aapxvi\ and the product and coproduct \aapxiv, \aapxv\ are
replaced by
$$(\aa\otimes \hh)(\bb\otimes \gg)=\aa(\hh_{(1)}\RRR \bb)\bxi(\hh_{(2)}\otimes 
\gg_{(1)})
\otimes \hh_{(3)}\gg_{(2)}\quad,\eqn\aapxxix$$
$$\Delta(\aa\otimes \hh)=\aa_{(1)}\bpsi(\hh_{(1)})^{(1)}\otimes 
\hh_{(2)}^{\ (1)}
\otimes \aa_{(2)}\bpsi(\hh_{(1)})^{(2)}\hh_{(2)}^{\ (2)}\otimes \hh_{(3)}
\quad.\eqn\aapxxx$$

It is convenient to have the explicit expression of \aapxxviii\ in the more 
simple cases.
For $\bpsi$ trivial it reads
$$
\Delta\bxi(\hh\otimes\gg)=
\bxi(\hh_{(1)}^{\ (1)}\otimes \gg_{(1)}^{\ (1)})\otimes\hh_{(1)}^{\ (2)}
(\hh_{(2)}\RRR\gg_{(1)}^{\ (2)})\bxi(\hh_{(3)}\otimes\gg_{(2)})
\quad.
\eqn\aapxxxi
$$
For $\RRR$ trivial, it gives
$$
\eqalign{
&
\Delta\bxi(\hh_{(1)}\otimes \gg_{(1)})\bpsi(\hh_{(2)}\gg_{(2)})=
\cr 
&
\mskip 50mu
\bpsi(\hh_{(1)})
[\bpsi(\gg_{(1)})^{(1)}\bxi(\hh_{(2)}^{\ (1)}\otimes \gg_{(2)}^{\ (1)})\otimes
\bpsi(\gg_{(1)})^{(2)}\hh_{(2)}^{\ (2)}\gg_{(2)}^{\ (2)}
\bxi(\hh_{(3)}\otimes \gg_{(3)})\quad.
\cr}
\eqn\aapxxxii
$$
For $\bpsi$ and $\RRR$ trivial, it reduces to
$$
\Delta\bxi(\hh\otimes \gg)=
\bxi(\hh_{(1)}^{\ (1)}\otimes \gg_{(1)}^{\ (1)})
\otimes \hh_{(1)}^{\ (2)}\gg_{(1)}^{\ (2)}\bxi(\hh_{(2)}\otimes \gg_{(2)})
\quad.
\eqn\aapxxxiii
$$

For $\bxi$ trivial $[\bxi(\hh\otimes \gg)=\epsilon(\hh)\epsilon(\gg)1_{{\AA}}]$ 
\aapxxii\ reduces
to \aapiii, \aapxxvii\ to \aapxii\ and \aapxxix\ to \aapxiv.
For $\bpsi$ trivial $[\bpsi(\hh)=1_{{\AA}}\otimes 1_{{\AA}}\epsilon(\hh)]$, 
\aapxxv\ reduces to \aapvi, \aapxxvi\ to \aapx\ 
and \aapxxx\ to \aapxv.
For $\beta(\hh)=\hh\otimes 1$ trivial, \aapxxx\ gives for the elements of $\HH$ 
with original primitive coproduct the cocycle extension expression
$$
\Delta(1\otimes\hh)=1\otimes\hh\otimes 1\otimes 1+
1 \otimes 1\otimes 1\otimes\hh+\bpsi(\hh)^{(1)}\otimes 1 
\otimes\bpsi(\hh)^{(2)}\otimes 1\quad,
\eqn\aapxxxx
$$
which in $\KK$ simply reads
$\Delta(\hh)=1\otimes\hh+\hh\otimes 1+\bpsi(\hh)$.
This was used for \hwrho\ [\insertii], \duali\ [\dualii] and \dualx\ [\dualxi].

\refout
\endpage
\end